\documentclass[onefignum,onetabnum]{siamart190516}
\usepackage{tikz}
\usetikzlibrary{arrows,decorations.pathmorphing,backgrounds,positioning,fit,petri}
\usetikzlibrary{calc,intersections,through,backgrounds}
\usepackage{subcaption}

\usepackage{booktabs}
\usepackage{lscape}

\usepackage{pstricks}
\usepackage{pst-node}
\usepackage{picinpar}
\usepackage{siunitx}
\usepackage{amssymb,amsmath,a4wide,mdwlist}
\usepackage{enumerate}
\usepackage{graphicx}
\usepackage{url}
\usepackage{rotating}

\usetikzlibrary{calc,intersections,through,backgrounds}

\newcommand{\Real}{\mathbb{R}}

\newtheorem{assum}{Assumption}[section]

\usepackage{lipsum}
\usepackage{amsfonts}
\usepackage{graphicx}
\usepackage{epstopdf}
\usepackage{algorithmic}
\ifpdf
  \DeclareGraphicsExtensions{.eps,.pdf,.png,.jpg}
\else
  \DeclareGraphicsExtensions{.eps}
\fi

\newsiamremark{remark}{Remark}
\newsiamremark{hypothesis}{Hypothesis}
\crefname{hypothesis}{Hypothesis}{Hypotheses}
\newsiamthm{claim}{Claim}
\newsiamremark{obs}{Observation}
\newsiamremark{example}{Example}

\headers{Sharp and fast bounds for CDT}{L.Consolini and M.Locatelli}

\title{Sharp and fast bounds for the Celis-Dennis-Tapia problem\thanks{Submitted to the editors DATE.}}

\author{Luca Consolini
\and Marco Locatelli\thanks{Dipartimento di Ingegneria e Architettura, Universit\`a di Parma
  (\email{luca.consolini,marco.locatelli@unipr.it}).}}

\usepackage{amsopn}

\title{Sharp and fast bounds for the Celis-Dennis-Tapia problem}

\date{}
\begin{document}
\maketitle
\begin{abstract}
In the Celis-Dennis-Tapia (CDT) problem a quadratic function is minimized over a region defined by two strictly convex quadratic constraints.
In this paper we re-derive a necessary and sufficient optimality condition for the exactness of the dual Lagrangian bound (equivalent to the Shor relaxation bound in this case).
Starting from such condition, we propose to strengthen the dual Lagrangian bound by adding one or two linear cuts to the Lagrangian relaxation. Such cuts are obtained from supporting
hyperplanes of one of the two constraints. Thus, they are redundant for the original problem but they are not for the Lagrangian relaxation.
The computational experiments show that the new bounds are effective and require limited computing times. In particular, one of the proposed bounds is able to solve
all but one of the 212 hard instances of the CDT problem presented in \cite{Burer13}.
\end{abstract}

\section{Introduction}
The Celis-Dennis-Tapia problem (CDT problem in what follows) is defined as follows:
\begin{equation}
\label{eq:cdt}
\begin{array}{ll}
p^\star=\min & {\bf x }^\top {\bf Q} {\bf x}+{\bf q}^\top {\bf x} \\ [6pt]
 & {\bf x}^\top {\bf x}\leq 1\\ [6pt]
& {\bf x}^\top {\bf A}{\bf x}+{\bf a}^\top {\bf x}\leq a_0,
\end{array}
\end{equation}
where ${\bf Q}, {\bf A}\in \mathbb{R}^{n\times n}$, ${\bf q},{\bf a}\in \mathbb{R}^n$, $a_0\in \mathbb{R}$,
while ${\bf A}$ is assumed to be positive definite. 
We will denote by 
$$
H=\{{\bf x}\in \mathbb{R}^n\ :\ {\bf x}^\top {\bf A}{\bf x}+{\bf a}^\top {\bf x}\leq a_0\},
$$ 
the ellipsoid defined by the second constraint, by $\partial H$ its border, and by $int(H)$ its interior.
The CDT problem was originally proposed in \cite{Celis85} and has 
attracted a lot of attention in the last two decades. For some special cases a convex reformulation is available. For instance, in \cite{Ye03} it is shown that a semidefinite reformulation is possible when no linear terms are present, i.e., when ${\bf q}={\bf a}={\bf 0}$. However, up to now no tractable convex reformulation of general CDT problems has been proposed in the literature. In spite of that, recently three different works \cite{Bienstock16,Consolini17,Sakaue16} independently proved that
the CDT problem is solvable in polynomial time. More precisely, in \cite{Consolini17,Sakaue16} polynomial solvability is proved by identifying  all KKT points through the solution of a bivariate
polynomial system with polynomials of degree at most $2n$. The two unknowns are the Lagrange multipliers of the two quadratic constraints. 
Instead, in \cite{Bienstock16} an approach based on the solution of a sequence of feasibility problems for systems of quadratic inequalities is proposed. The systems are solved 
by a polynomial-time algorithm based on Barvinok's construction  \cite{Barvinok93}. 
Though polynomial, all these approaches are computationally demanding since the degree of the polynomial is quite large.
Conditions guaranteeing that the classical Shor SDP relaxation or, equivalently in this case, the dual Lagrangian bound is exact, are discussed in \cite{Ai09,Beck06}. In particular, in \cite{Ai09} a necessary and sufficient condition is presented. It is shown that lack of exactness is related to the existence of KKT points with the same Lagrange multipliers but two distinct
primal solutions, both active at one of the two constraints but one violating and the other one fulfilling the other constraint.
In \cite{Bomze15} necessary and sufficient conditions for local and global optimality are discussed based on copositivity.
In \cite{Bomze18} an exactness condition is given for a copositive relaxation, also for the case with additional linear constraints.
A trajectory following method to solve the CDT problem has been discussed in \cite{Ye03}, while different branch-and-bound solvers are tested in \cite{Montahner18}.
\newline\newline\noindent
Recently, different papers proposed valid bounds for the CDT problem.
In \cite{Burer13} the Shor relaxation bound is strengthened by adding all RLT constraints obtained by supporting hyperplanes of the two ellipsoids. By fixing the supporting hyperplane for one ellipsoid, the RLT constraints obtained with all the supporting hyperplanes of the other can be condensed into a single SOC-RLT constraint. Varying the supporting hyperplane of the first ellipsoid gives rise to an infinite number of SOC-RLT constraints which, however, can be separated in polynomial time. The addition of these constraints does not allow to close the duality gap but it is computationally shown that many instances which are not solved via the SDP bound, are solved with the addition of these SOC-RLT cuts. The authors generate 1000 random test instances for each $n=5,10,20$,
following a procedure described in \cite{Martinez94} to generate trust-region problems with one local and nonglobal minimizer. The proposed bound based on SOC-RLT cuts allows to solve most instances except for 212 (38 for $n=5$, 70 for $n=70$, and 104 for $n=20$). Such unsolved instances are considered as hard ones in subsequent works.
In \cite{Yang16} lifted-RLT cuts are introduced and it is shown that the new constraints allow to derive an exact bound for $n=2$ but also to improve the bounds of \cite{Burer13} over the hard instances for $n>2$.
In \cite{Yuan17} it is proved that the duality gap can be reduced to 0 by solving two subproblems with SOC constraints when the second constraint is the product of two linear functions.
% A model with two SOC constraints is also presented which, in some cases, allows to close the duality gap......CONTROLLARE\newline
In \cite{Anstreicher17} KSOC cuts are introduced. These are Kronecker product constraints which generalize both the classical RLT constraints obtained from two linear inequality constraints, and the SOC-RLT constraints obtained 
from one linear inequality constraint and a SOC constraint. Further hard instances from \cite{Burer13} are solved with the addition of these cuts. 
%\newline\newline\noindent
%In \cite{Yang16} all the constraints which are necessary to close the gap are introduced for the case when $n=2$ but not for larger dimensions, although new constraints allow to reduce the gap.
%\newline\newline\noindent
%In \cite{Yuan17} it is proved that the duality gap can be reduced to 0 by solving two subproblems with SOC constraints when the second constraint is the product of two linear functions. A model with two SOC constraints is also presented which, in some cases, allows to close the duality gap.
\newline\newline\noindent
In this paper we investigate ways to strengthen the dual Lagrangian bound through the addition of one or two linear cuts. In particular,
the paper is structured as follows. 
In Section \ref{sec_gen_prop} we derive some theoretical results for a class of problems with two constraints which includes the CDT problem as a special case. 
We develop a bisection technique to solve the dual Lagrangian relaxation for such class of problems. 
In the following sections we apply the results of Section \ref{sec_gen_prop} to the CDT problem.
In particular, in Section \ref{sec:dual} we introduce some results through which it will be possible to re-derive the necessary and sufficient exactness condition discussed in \cite{Ai09}.
In Section \ref{sec:boundimp} we discuss how to improve the dual bound for the CDT problem by the addition of a linear cut. Next, in Sections \ref{sec:optimal}-\ref{sec:two} we discuss techniques to further improve the bound. More precisely, in Section \ref{sec:optimal} we still present a bound based on the addition of a linear cut but we develop a technique to locally adjust a given linear cut, while in Section \ref{sec:two} we consider a bound based on the addition of two linear cuts.
Finally, in Section \ref{sec:comp} we present some computational experiments which show that the newly proposed bounds, in particular those based on two linear cuts, are both computationally cheap and effective. In particular, one of the bounds will be able
 to solve all but one of the hard instances from \cite{Burer13}. We also investigate which are the most challenging instances for the proposed bounds and, as we will see, the difficulties are related to the existence of multiple solutions
 of Lagrangian relaxations.

\section{Lower bounds obtained from the Lagrangian relaxation}
\label{sec_gen_prop}

The CDT problem~\eqref{eq:cdt} is a specific instance of the
following, more general, one:

\begin{equation}
\label{eq:gen}
\begin{array}{ll}
p^\star=\min_{{\bf x} \in \mathbb{R}^n} & f({\bf x}) \\ [6pt]
 & g({\bf x})\leq 0,\\ [6pt]
& h({\bf x}) \leq 0.
\end{array}
\end{equation}
In this section, we discuss a class of lower bounds on the solution of
problem~\eqref{eq:gen} that can be obtained from its Lagrangian
relaxation. In the next sections, we will apply these bounds to the
specific case of the CDT problem~\eqref{eq:cdt}. Throughout this and the following sections, we make the following assumptions.
\begin{assum}
  \label{assum_p}
  In Problem~\eqref{eq:gen}
\begin{description}

\item[a)] $g,h$ are continuous;

\item[b)] the set $\{{\bf x} \in  \mathbb{R}^n: g({\bf x}) \leq 0\}$ is bounded;

\item[c)] the problem is strictly feasible, that is
\begin{equation}
\label{eq:ellag}
h_0=\min_{{\bf x}\ :g({\bf x})
  \leq 0 } h({\bf x}) < 0;
\end{equation}

\item[d)] the solution set of Problem~\eqref{eq:gen} without the last
  constraint, that is
\[
  \begin{array}{ll}
\bar P=\arg \min_{{\bf x} \in \mathbb{R}^n} & f({\bf x}) \\ [6pt]
 & g({\bf x})\leq 0,
\end{array}
\]
is such that $(\forall {\bf x} \in \bar P)\; h({\bf x}) >0$. In other words,
constraint $h({\bf x}) \leq 0$ is active on the solution set of~\eqref{eq:gen}.
\end{description}
\end{assum}

Note that if the last condition in Assumption~\ref{assum_p} is
violated, we can find
the solution of Problem~\eqref{eq:gen} by removing the last
constraint and the relaxation discussed in this section is useless.
Now, let $G=\{{\bf x} \in \Real^n : g({\bf x}) \leq 0\}$ and $H=\{{\bf x} \in \Real^n : h({\bf x}) \leq 0\}$.
Let $X\supset H$ be a closed subset of $\Real^n$ and for $\lambda \in \Real$,
with $\lambda \geq 0$, define the Lagrangian relaxation
\begin{equation}
\label{eq:paramcdtg}
p_X(\lambda)=\min_{{\bf x}\in X \cap G}  f({\bf x}) + \lambda h({\bf x}),
\end{equation}
and the corresponding solution set
\[
P_X(\lambda)=\arg\min_{{\bf x}\in X \cap G}  f({\bf x}) + \lambda h({\bf x}).
\]

Note that $P_X(\lambda)$ is compact, since $G\cap X$ is nonempty (in view of part c) of Assumption~\ref{assum_p}) and compact (in view of the compactness of $G$ which follows from parts a) and b) of  Assumption~\ref{assum_p}), and
$f+ \lambda h$ is continuous. Due to well-known properties of the Lagrangian
relaxation, we have that function $p_X$ is such that $(\forall \lambda \geq 0)\; p_X(\lambda) \leq p^*$,
and is concave (it is the pointwise minimum of a set of functions linear in $\lambda$).
The best bound that can be obtained as the solution
of~\eqref{eq:paramcdtg} is given by
\begin{equation}
  \label{eqn_dual}
\bar{p}_X=\max_{\lambda \geq 0} p_X(\lambda),
\end{equation}
and corresponds to the solution of the dual Lagrangian problem.
Note that function $p_X$ depends on the choice of set $X$. 

Now, we recall that the supergradient of a function $q: \Real \to \Real$ at $x
\in \Real$, is defined as
\[
\partial q(x)= \{z \in \Real: (\forall y \in \Real) \, q(y)-q(x) \leq z(y-x)\}.
\]
Since $p_X$ is concave, for any $\lambda \in \Real$, the supergradient
$\partial p_X(\lambda)$ is non-empty.

For $A \subset \mathbb{R}^n$ define the following subset of $\mathbb{R}$
\[
h(A)=\{h({\bf x}), {\bf x} \in A\}.
\]
For $X \subset \mathbb{R}^n$, define a (set-valued) function $Q_X: \mathbb{R}^+ \to \mathcal{P}(\mathbb{R})$
\begin{equation}
  \label{eq_Q}
Q_X(\lambda)=h(P_{X}(\lambda))
\end{equation}
($\mathbb{R}^+$ denotes the set of nonnegative reals and $\mathcal{P}(\mathbb{R})$ is the power set of the set of real numbers).
Also set $h^{\min}_X(\lambda)=\min Q_X(\lambda)$,
$h^{\max}_X(\lambda)=\max Q_X(\lambda)$.
%Further, set $\lambda_X=\min\{\lambda \geq 0: 0 \in \partial p_X(\lambda)\}$.
%$h^{\min}_X= h^{\min}_X(\lambda_X)$,
%$h^{\max}_X= h^{\max}_X(\lambda_X)$.
The following proposition shows that function $Q_X$ is monotone non-increasing (see Definition~3.5.1
of~\cite{aubin2009set}) and upper semicontinuous (see Definition~1.4.1 of~\cite{aubin2009set}).
These two properties will play an important role in the computation of
a lower bound for Problem~\eqref{eq:gen}. Moreover, this proposition characterizes the supergradient of $p_X$
at each $\lambda\geq 0$. 
\begin{proposition}
  \label{prop_for_Q}
  For any $X \subset \Real^n$

  i) $Q_X$ is monotone not-increasing, that is
if $\lambda_1 \geq \lambda_2$,
  $y_1 \in Q_X(\lambda_1)$, $y_2 \in Q_X(\lambda_2)$, then $y_1 \leq
  y_2$.

ii) $Q_X$ is upper semicontinuous, that is, if $Q_X(\lambda) \subset
U$,
where $U$ is an open subset of $\Real$, then there exists a neighborhood $V$ of
$\lambda$ such that $(\forall z \in V) \, Q_X(z) \subset U$.

iii) $\partial p_X(\lambda)=[\min
Q_X(\lambda), \max Q_X(\lambda)]$.

  \end{proposition}
\begin{proof}
  i) Let ${\bf x}_1,{\bf x}_2 \in \Real^n$ be such that $y_1=h({\bf x}_1)$ and
  $y_2=h({\bf x}_2)$, then
  $f({\bf x}_1) + \lambda_1 h({\bf x}_1) \leq f({\bf x}_2) + \lambda_1 h ({\bf x}_2)$ and
  $f({\bf x}_2) + \lambda_2 h({\bf x}_2) \leq f({\bf x}_1) + \lambda_2 h({\bf x}_1)$.
  By adding up the previous inequalities, it follows that
  $(\lambda_1-\lambda_2) (h({\bf x}_1)-h({\bf x}_2)) \leq 0$, which implies the thesis.

ii) It is a consequence of Berge's Maximum Theorem
(see~\cite{Berge1963}). In particular, we consider the slightly
different formulation presented as the corollary to Theorem~3 on page
30 of~\cite{Hil74}. Namely, since function $(G \cap X) \times \Real \to \Real$, $({\bf x},\lambda) \leadsto f({\bf x}) + \lambda
h({\bf x})$ is continuous, set valued function $P_X$ is upper
semicontinuous. Hence, also $Q_X$ is upper semicontinuous, since it is
obtained as the composition of $P_X$ with $h$, which is continuous
(see Theorem~1' on page~113 of~\cite{Berge1963}).

iii) It is a consequence of Theorem~4.4.2 in~\cite{ConAl}, being $G$ compact.

\end{proof}
The next proposition characterizes the optimal solution of the dual Lagrangian problem (\ref{eqn_dual}).
\begin{proposition}
Under Assumption~\ref{assum_p} the optimal value of the dual Lagrangian problem (\ref{eqn_dual}) is attained at
$\lambda_X>0$ such that $0 \in \partial p_X(\lambda_X)$.
%$$
%\lambda_X=\min\{\lambda >0: 0 \in \partial p_X(\lambda)\}.
%$$
\end{proposition}
\begin{proof}
We first prove that the value $\bar{p}_X$ is attained. Let ${\bf x}_0$ be an optimal solution of problem (\ref{eq:ellag}).
Since ${\bf x}_0\in H$ (more precisely, it belongs to the interior of $H$) and recalling that $H\subset X$, it holds that ${\bf x}_0\in G\cap X$. 
Then, for each $\lambda$
$$
p_X(\lambda)\leq f({\bf x}_0)+\lambda h({\bf x}_0)\rightarrow -\infty \ \ \ \mbox{as}\ \ \ \lambda\rightarrow +\infty
$$
(recall that $h({\bf x}_0)=h_0<0$ in view of part c) of Assumption~\ref{assum_p}). Thus, the maximum value of $p_X(\lambda)$ is attained at some $\lambda_X\geq 0$. But in view of part d) of Assumption~\ref{assum_p}, we have that  $h^{\min}_X(0)>0$. Thus, function $p_X$ is increasing at $\lambda=0$ and, consequently, we must have $\lambda_X>0$. Moreover, by the optimality condition of nondifferentiable concave functions,  $0\in \partial p_X(\lambda_X)$ must hold.
\end{proof}
The following property shows that it is always possible to find a
sufficiently high
value of $\lambda$ such that
$P_X(\lambda) \subset H$, that is, the elements of $P_X(\lambda)$ are
feasible solutions of Problem~\eqref{eq:gen}.

\begin{lemma}
  \label{lem:2g}
If 
\begin{equation}
\label{eq:condlambda_g}
\lambda\geq \hat{\lambda}=\frac{\max_{{\bf x}\in G\cap X}  f({\bf x})-
\min_{{\bf x} \in G \cap X}  f({\bf x})}
{|h_0|},
\end{equation}
where $h_0$ is defined in (\ref{eq:ellag}), then
$P_{X} (\lambda) \subset H$.
\end{lemma}

\begin{proof}
  By contradiction, assume that there exists ${\bf x} \in P_X(\lambda)$
  such that $h({\bf x}) >0$ and let ${\bf x}_0 \in G \cap H$ be
  such that $h({\bf x}_0)=h_0<0$,
  then
  $f({\bf x})+ \lambda h({\bf x}) \leq f({\bf x}_0)+ \lambda h({\bf x}_0)$. Since $h({\bf x})>0$, it
  follows that $\lambda \leq \frac{f({\bf x}_0)-f({\bf x})}{|h({\bf x}_0)|+h({\bf x})}
  < \frac{\max_{{\bf x}\in G\cap X}  f({\bf x})-
\min_{{\bf x} \in G \cap X}  f({\bf x})}{|h({\bf x}_0)|}$ which contradicts the assumption
on $\lambda$.
  \end{proof}

%  If $Q_X(\lambda)$ contains negative elements, then
%the corresponding elements of $P_X(\lambda)$ are feasible solutions for~\eqref{eq:gen}:
   
 % \begin{proposition}
  %  \label{prop_low_bound}
   % $P_X(\lambda)$ contains a
%feasible solution for Problem~\eqref{eq:gen} if and only if $\min
%Q_X(\lambda) \leq 0$.
 %   \end{proposition}

%    \begin{proof}
%$P_X(\lambda)$ contains a feasibile solution for~\eqref{eq:gen} if and
%only if there exists $x \in P_X(\lambda)$ such that $h(x) \leq 0$,
%which is true if and only if $\min Q_X(\lambda)\leq 0$.
% \end{proof}

 The following proposition shows that if $0 \in Q_X(\lambda)$, then
$p_X(\lambda)$ is equal to the optimal value of Problem~\eqref{eq:gen}. 
 
\begin{proposition}
  \label{prop:noexact_g}
Under Assumption~\ref{assum_p}
the following statements are equivalent:
  
i) $0 \in Q_X(\lambda)$,

ii) $p^* = p_{X} (\lambda)$ and there exists $\bar {\bf x} \in \arg \min_{{\bf x} \in
  G \cap H} f({\bf x})$ such that $h(\bar {\bf x})=0$.
\end{proposition}

\begin{proof}
 
i) $\Rightarrow$ ii).  Let $\bar {\bf x}$ be such that $h(\bar {\bf x}) \in Q_X(\lambda)$ and
$h(\bar {\bf x})=0$. Let ${\bf x}^*$ be a solution of~\eqref{eq:gen}.
 Then, $p_X(\lambda)= f(\bar {\bf x}) + \lambda h(\bar {\bf x}) =f(\bar {\bf x})
 \leq  f({\bf x}^*) + \lambda h({\bf x}^*) \leq f({\bf x}^*)$, hence
 $p_X(\lambda) \leq p^*$. Moreover, $p_X(\lambda)=f(\bar {\bf x}) \geq \min_{{\bf x} \in
   G \cap H} f({\bf x}) = p^*$. 

ii) $\Rightarrow$ i). Assume that $p_{X} (\lambda)=p^*$ and
let ${\bf x} \in P_X(\lambda)$.
Then, by ii), $f({\bf x}) + \lambda h({\bf x})  = f(\bar {\bf x})= f(\bar {\bf x}) + \lambda
h(\bar {\bf x})$.
It follows that $\bar{\bf x} \in P_X(\lambda)$ and $Q_X(\lambda) \ni h(\bar{\bf x})=0$.
\end{proof}

\begin{remark}
  If $0 \in Q_X(\lambda)$, by point iii) of Proposition~\ref{prop_for_Q}, $\partial
  p_X(\lambda) \ni 0$, so that $\lambda$ corresponds to a maximizer of
  the dual Lagrangian. Note that equation $0 \in Q_X(\lambda)$
  always admits a solution if $Q_X$ is continuous. However,
  in the general case, $Q_X$ is only upper semicontinuous.
In this case, the value of $\lambda$ for which $\partial p_X(\lambda)
\ni 0$  may not satisfy $0 \in Q_X(\lambda)$. Thus, the
optimal value of the dual Lagrangian~\eqref{eqn_dual} is not equal to the optimal value of~\eqref{eq:gen} but it represents a lower bound of it.
\end{remark}
%For ease of notation, set
%$p_X=p_X(\lambda_X)$.
In order to evaluate a numerical solution algorithm, we define the
following weak solution of~\eqref{eq:gen}.
\begin{definition}
  ${\bf x}$ is an $\eta$-solution of~\eqref{eq:gen} if
 ${\bf x} \in G \cap H$ and $f({\bf x})-p^* \leq \eta$.
\end{definition}
The following proposition presents a bound on the error committed on
the estimation of $p^*$.
%\begin{proposition}
 % \label{prop_est}
%For any $\lambda \geq 0$, there exists $x \in P_X(\lambda)$ such that
%  $f(x)-p^* \leq \lambda \, \min\{|z|: z \in Q_X(\lambda)\}$.
%\end{proposition}
%\begin{proof}
%Let $x \in P_X(\lambda)$ be such that $|h(x)|=\min\{|z|: z \in Q_X(\lambda)\}$.
%Then $f(x) + \lambda h(x) \leq f(x^*)+ \lambda h(x^*)
%\leq f(x^*)$,
%from which $f(x)-f(x^*) \leq \lambda |h(x)|$. 
%\end{proof}
\begin{proposition}
  \label{prop_est}
  For any $\lambda \geq 0$ such that $ P_X(\lambda)\cap H\neq \emptyset$, and for any ${\bf x} \in P_X(\lambda)\cap H$,
it holds that 
  $f({\bf x})-p^* \leq \lambda |h({\bf x})|$, i.e., ${\bf x}$ is an $\eta$-solution of problem of~\eqref{eq:gen} with $\eta=\lambda |h({\bf x})|$.
\end{proposition}
\begin{proof}
Since ${\bf x} \in P_X(\lambda)$ and observing that ${\bf x}^*\in G\cap X$ for any $X\supset H$,
$f({\bf x}) + \lambda h({\bf x}) \leq f({\bf x}^*)+ \lambda h({\bf x}^*)
\leq f({\bf x}^*)$,
from which $f({\bf x})-f({\bf x}^*) \leq \lambda |h({\bf x})|$. 
\end{proof}
Now we introduce Algorithm~\ref{algo:1} which is based on a binary search through different $\lambda$ values and is able to return the solution of the dual Lagrangian problem, i.e., the maximum of function $p_X(\lambda)$ and, in some cases, even the solution of problem (\ref{eq:gen}). %The algorithm is described for a generic region $X\supseteq {\cal E}$, while later on we will comment the special case $X=\mathbb{R}^n$. 
%Algorithm~\ref{algo:1}, based on a binary
%search through different $\lambda$ values and is able to return the
%solution of the dual Lagrangian problem, i.e., the maximum of function
%$p_X(\lambda)$ and, in some cases, even the solution of
%problem~\eqref{eq:gen}. The algorithm is described for a generic region $X\supseteq {\cal E}$, while later on we will comment the special case $X=\mathbb{R}^n$. 
The algorithm also returns a point ${\bf z}_1(\lambda^{\max})\in H$  and  (possibly) a point ${\bf z}_2(\lambda^{\min})\not\in H$. Note that according to Proposition \ref{prop_est}, point ${\bf z}_1(\lambda^{\max})$ is an $\eta$-solution of problem of~\eqref{eq:gen} with $\eta=\lambda |h({\bf z}_1(\lambda^{\max}))|$.
%the set $Z_X=\{x \in X: h(x)=\max Q_X(\lambda_X)\}$ 
%that will be used later on.

\begin{algorithm}
\caption{Binary search algorithm for the solution of the dual Lagrangian problem for (\ref{eq:cdt}).}
\label{algo:1}
%\AlgData{${\bf Q}, {\bf A}\in \mathbb{R}^{n\times n}$, ${\bf q},{\bf a}\in \mathbb{R}^n$, $a_0\in \mathbb{R}$, $\varepsilon>0$\;}
%\AlgResult{$\boldsymbol{\gamma}$}
\vspace*{0.25cm}
%\hspace*{\algorithmicindent} 
\textbf{DualLagrangian}($X$, $\lambda_{\tt init}$)
\vspace*{0.25cm}
\begin{algorithmic}
\STATE{Set $\lambda_{\min}=0$, $\lambda^{\max}=\lambda^{\tt init}$}
\WHILE{$\lambda^{\max}-\lambda^{\min}>\varepsilon$}
\STATE{Set $\lambda=(\lambda^{\max}+\lambda^{\min})/2$}
\STATE{Solve problem (\ref{eq:paramcdt}) and let $P_X(\lambda)$ be its set of optimal solutions}
\STATE{Compute the set $Q_X(\lambda)$ and the values $h_X^{\min}(\lambda),h_X^{\max}(\lambda)$}
\IF{$h_X^{\min}(\lambda)>0$}
\STATE{Set $\lambda_{\min}=\lambda$}
\ELSIF{$h_X^{\max}(\lambda)<0$}
\STATE{Set $\lambda^{\max}=\lambda$}
\ELSE
\STATE{Set $\lambda^{\max}=\lambda^{\min}=\lambda$}
\ENDIF
\ENDWHILE
\STATE{Set $Lb=p_X(\lambda^{\max})$, let ${\bf z}_1(\lambda^{\max})$ be some point in $P_X(\lambda^{\max})\cap H$ and ${\bf z}_2(\lambda^{\min})$ be some point (if any) in $P_X(\lambda^{\min})\setminus H$}
%\RETURN $Lb,\lambda^{\min},\lambda^{\max},{\bf z}_1(\lambda^{\max}), {\bf z}_2(\lambda^{\min}),h_X^{\min},h_X^{\max}$
\RETURN $Lb,\lambda^{\max},{\bf z}_1(\lambda^{\max}), {\bf z}_2(\lambda^{\min})$
\end{algorithmic}
\end{algorithm}
The algorithm starts with an initial interval of $\lambda$ values $\left[\lambda^{\min},\lambda^{\max}\right]=[0,\lambda^{\tt init}]$, where $\lambda^{\tt init}$ is a suitably large value and can be set equal to $\hat{\lambda}$ as defined in Lemma~\ref{lem:2g}.
%defined in (\ref{eq:condlambda}). 
At each iteration the algorithm halves such interval 
by evaluating the set $Q_X^{\lambda}$ at $\lambda=(\lambda^{\max}+\lambda^{\min})/2$. Then, the algorithm sets: $\lambda^{\min}=\lambda$, if $h_X^{\min}(\lambda)>0$; $\lambda^{\max}=\lambda$ if
$h_X^{\max}(\lambda)<0$. Instead, if $0\in \partial
p_X(\lambda)=[h_X^{\min}(\lambda),h_X^{\max}(\lambda)]$, the
algorithm sets 
$\lambda^{\max}=\lambda^{\min}=\lambda$ and exits the loop.

The following proposition characterizes Algorithm~\ref{algo:1}.

\begin{proposition}
\label{prop_alg}
  \mbox{  }
  
i)   Algorithm~\ref{algo:1} terminates in a finite number of iterations,

ii) at each iteration $\lambda^{\min} \leq \lambda_X \leq \lambda^{\max}$,

iii) at termination $|\lambda^{\max}-\lambda_X| \leq \epsilon$,

iv) at each iteration, if $\lambda^{\min} < \lambda_X  < \lambda^{\max}$, then $[h_X^{\max}(\lambda^{\max}), h_X^{\min}(\lambda^{\min})] \supset
\partial p_X(\lambda_X)$,

v) point ${\bf z}_1(\lambda^{\max})\in P_X(\lambda^{\max})\cap H$ is an $\eta$-solution of~\eqref{eq:gen}
with
$\eta=\lambda_{\max}|h({\bf z}_1(\lambda^{\max}))|$.
%$ \, \min\{|h_X^{\min}(\lambda^{\max})|,|h_X^{\max}(\lambda^{\min})|\}$.
\end{proposition}

\begin{proof}
i)  At each iteration the length of the interval
  $\left[\lambda^{\min},\lambda^{\max}\right]$ is halved. Hence, in
  a sufficient large number of iterations, the termination condition
  of the main loop will be satisfied.

ii) At the beginning of
the algorithm we have that $\lambda^{\min} \leq \lambda_X \leq
\lambda^{\max}$.
Every time $\lambda^{\min}$ is updated, we set
$\lambda^{\min}=\lambda$ if condition
$h_X^{\min}(\lambda)>0$ holds. Since $h_X^{\min}(\lambda_X) \leq  0$, 
by the monotonicity of function $h_X^{\min}$, which is a consequence of the monotonicity of function $Q_X$,
condition $\lambda^{\min} \leq
\lambda_X$ is maintained. The same reasoning can be used to prove that
$\lambda^{\max} \geq \lambda_X$.

iii) It is a consequence of ii) and the termination condition.  

iv) $\partial p_X(\lambda_X)=
[h_X^{\min}(\lambda_X),h_X^{\max}(\lambda_X)] \subset
[h_X^{\max}(\lambda^{\max}), h_X^{\min}(\lambda^{\min})]$,
  due to point ii) and the monotonicity of functions
  $h_X^{\max}$ and $h_X^{\min}$, which is a consequence of the monotonicity of function $Q_X$.

v) It is a consequence of Proposition~\ref{prop_est}. % and point iv).
  
\end{proof}

%As a consequence of the previous property (iv), if the algorithm
%terminates with either $h_{\min}^{\lambda_{\max}}(X)=0$ or
%$h_{\max}^{\lambda_{\min}}(X)=0$, then the solution is exact.
%Note that, since $\lambda_{\min}$ and $\lambda_{\max}$ are continuous
%functions, the estimate of the superdifferential is arbitrarily
%precise (DA SISTEMARE).

%If we remove the exit condition $\lambda^X_{\max}-\lambda^X_{\min}
%\leq \varepsilon$, then either the algorithm stops after a finite number of iterations, as soon as the $exit$ flag is set equal to ${\tt true}$, or it runs for an infinite number of iterations. In the latter case we set $\bar{\lambda}^X$ equal to the common limit of $\lambda^X_{\min},\lambda^X_{\max}$, i.e.,
%$\lambda^X_{\min},\lambda^X_{\max}\rightarrow \bar{\lambda}^X$. For such limit value we have that
%$$
%h_{\min}^{\bar{\lambda}^X}(X)=\lim h_{\max}^{\lambda_{\max}^X}(X),\ \ \ h_{\max}^{\bar{\lambda}^X}(X)=\lim h_{\min}^{\lambda_{\min}^X}(X),
%$$
%where $0\in \left[h_{\min}^{\bar{\lambda}^X}(X),h_{\max}^{\bar{\lambda}^X}(X)\right]$, so that $\bar{\lambda}^X$ is the maximizer of $p_X$ and $p_X(\bar{\lambda}^X)$ is the dual Lagrangian bound. 

The following property is a direct consequence of the upper
semicontinuity of $Q_X$.
\begin{proposition}
  \label{prop:improveg}
  Let $X \supset H$ be such
  that $\sup Q_X(\lambda) <0$, then
  there exists a neighborhood $U$ of $\lambda$ such
  that $(\forall \eta \in U)\; \max Q_X(\eta) <0$.
\end{proposition}

As a consequence of the previous proposition, it is possible to improve the lower
bound on Problem~\eqref{eq:gen}, obtained as the solution
of~\eqref{eq:paramcdtg}, by 
replacing set $X$ with a different set $Y\supset H$ fulfilling a given condition. %such that $P_Y(\lambda_X)\setminus H=\emptyset$.
%restricting set $X$ to a smaller set
%$Y$. In particular, we want to remove from $X$ 
%the elements of $P_X(\lambda_X)$ at which function $h$ is nonnegative.

\begin{proposition}
\label{prop:boundimprove}
Let $Y \supset H$
 be such that $\max Q_{Y}(\lambda_{X})\leq 0$ or, equivalently, $P_Y(\lambda_X)\setminus H=\emptyset$, and assume that $\bar{p}_X=p_{X}(\lambda_X) < p^*$.
 Then
$\bar{p}_Y=p_{Y}(\lambda_Y) > \bar{p}_{X}$.
%In particular, if $Y\subset X$, then the result holds if $Y\cap(P_X(\lambda_X)\setminus H)=\emptyset$.
 \end{proposition}

\begin{proof}
If $\max Q_{Y}(\lambda_{X})=0$, then $0\in Q_{Y}(\lambda_{X})$ and by Proposition \ref{prop:noexact_g}
$\bar{p}_Y=p^*>\bar{p}_X$. Thus, we only consider the case $\max Q_{Y}(\lambda_X)<0$.
In such case, by Proposition~\ref{prop:improveg}, $\lambda_Y< \lambda_{X}$.
Since $0 \in [h_{Y}^{\min}(\lambda_Y), h_{Y}^{\max}(\lambda_Y)]$,
there exists
${\bf y} \in Q_Y(\lambda_Y)$ such that $h({\bf y}) \leq 0$.
Note that $\bar{p}_{Y}=f({\bf y})+ \lambda_Y h({\bf y})$.
If $h({\bf y})=0$, then, by Proposition~\ref{prop:noexact_g}, $p_Y(\lambda_Y)=p^*$, so
that the thesis is satisfied in view of $\bar{p}_X<p^*$.
Otherwise,  if $h({\bf y})<0$, let ${\bf x} \in \Real^n$ be such that $\bar{p}_{X}=f({\bf x})+ \lambda_X h({\bf x})$.
 Then $\bar{p}_{X}=f({\bf x}) + \lambda_{X}  h({\bf x}) \leq
  f({\bf y})+ \lambda_{X} h({\bf y}) <   f({\bf y})+ 
  \lambda_Y h({\bf y})=\bar{p}_{Y}$, where we used the facts that $h({\bf y})<0$
%, being $\max Q_{Y}(\lambda_{X})<0$, 
and that $\lambda_Y < \lambda_X$.
%Now, since $h_X^{\min}(\lambda_X)<0$ we have that $P_X(\lambda_X)\cap H\neq \emptyset$ and, consequently, since $Y\supset H$,  also $Y\cap P_X(\lambda_X)\neq \emptyset$.
%If $Y\subset X$, then $P_Y(\lambda_X)=Y\cap P_X(\lambda_X)$. Moreover, if 
%$Y\cap(P_X(\lambda_X)\setminus H)=\emptyset$, then the condition $\max Q_{Y}(\lambda_X)<0$ is satisfied.
\end{proof}
The following proposition deals with the special case of the previous result when $Y\subset X$.
\begin{proposition}
\label{prop:boundimprove1}
Let $X\supset Y \supset H$
 be such that 
$Y\cap(P_X(\lambda_X)\setminus H)=\emptyset$, and assume that $\bar{p}_X=p_{X}(\lambda_X) < p^*$.
 Then
$\bar{p}_Y=p_{Y}(\lambda_Y) > \bar{p}_{X}$.
 \end{proposition}
\begin{proof}
Since $h_X^{\min}(\lambda_X)<0$ we have that $P_X(\lambda_X)\cap H\neq \emptyset$ and, consequently, since $Y\supset H$,  also $Y\cap P_X(\lambda_X)\neq \emptyset$.
Then, $Y\subset X$ implies $P_Y(\lambda_X)=Y\cap P_X(\lambda_X)$. Moreover, if 
$Y\cap(P_X(\lambda_X)\setminus H)=\emptyset$, then the condition $\max Q_{Y}(\lambda_X)\leq 0$ is satisfied and the result follows from Proposition \ref{prop:boundimprove}.
\end{proof}
Stated in another way, the previous propositions show that, in case the lower bound $\bar{p}_X$ is not exact, we are able to improve (increase) it, if we are able to replace set $X$ with a new set $Y$ which
cuts away all members of $P_X(\lambda_X)$ outside $H$.
\begin{remark}
\label{rem:diffp}
Up to now we have not discussed 
the difficulty of computing the values of function $p_X$ or, equivalently, the difficulty of solving problem (\ref{eq:paramcdtg}). Such difficulty is strictly related to the specific problem (i.e., to the specific functions $f,g,h$), and also to the specific set $X$.
In the next sections we apply the general theory developed in this section to the CDT problem.
We show that for suitably defined sets $X$ (defined by one or two linear cuts), the computation of function $p_X$ can be done in an efficient way, and, moreover, the corresponding lower bounds $\bar{p}_X$ improve the standard dual Lagrangian bound, corresponding to the case $X=\mathbb{R}^n$. 
\end{remark}
\begin{remark}
In principle one could also define a cutting algorithm where a sequence of sets $\{X_k\}$ is generated such that: i) $X_k\supset X_{k+1} \supset H$ for all $k$; ii) $X_{k+1}\cap (P_{X_k}(\lambda_{X_k})\setminus H)=\emptyset$; iii) $\cap_{k=1}^\infty X_k =H$. The corresponding sequence of lower bounds $\{\bar{p}_{X_k}\}$ is strictly increasing in view of Proposition \ref{prop:boundimprove1}, and converges to $p^*$. However, the difficulty related to such an algorithm is that forcing ii) may not be trivial and, moreover, as already commented in Remark \ref{rem:diffp}, computing $p_{X_k}$ may be computationally demanding. 
\end{remark}

The following Algorithm~\ref{algo:3bis}, in principle, is able
to always find an approximate solution of~\eqref{eq:gen}.
The algorithm is based on an iterative reduction of set $X$, in
order to eliminate its elements in which function $h$ is positive.
In practice, Algorithm~\ref{algo:3bis} could be unimplementable. Indeed,
it may require a large number of cuts on set $X$ and each added cut may
increase the complexity of the optimization problem that we need to
solve to evaluate \textbf{DualLagrangian}. In
Section~\ref{sec:optimal}, we will see that, to refine the lower bound
on the solution of the CDT problem, it is computationally more
convenient to adjust existing cuts instead of adding new ones.
We stress that we will not actually use Algorithm~\ref{algo:3bis} for the
solution of the CDT problem. We present this algorithm just as a
theoretical contribution. 

\begin{algorithm}
\caption{Bound improvement through redefinition of set $X$.}
\label{algo:3bis}
\vspace{0.25cm}
%\hspace*{\algorithmicindent} 
\begin{algorithmic}[1]
\STATE{Set $X=\Real^n$}

\STATE{Set $\lambda^{\max}=\lambda^{\tt init}$}

\REPEAT
\STATE{Let $[Lb,\lambda^{\min},\lambda^{\max},{\bf z}_1(\lambda^{\max}), {\bf z}_2(\lambda^{\min}),h_X^{\min},h_X^{\max}]= \textbf{DualLagrangian}(X,\lambda^{\tt init})$} 
\STATE{Set $Z=\{{\bf x} \in P_X(\lambda^{\min}): h({\bf x})>0\}$}
\STATE{Redefine $X=Y$, where $Y$ is such that $X \supset Y \supset H$ and $Z
  \cap Y = \emptyset$. }
\UNTIL{$\min\{h_X^{\max}(\lambda^{\min}),-h_X^{\min}(\lambda^{\max})\}
  \lambda^{\max} \leq\eta$}
\RETURN{$\bar{\bf x} \in P_X(\lambda_X^{min}) \cup P_X(\lambda_X^{max})$
  with $|h(\bar{\bf x})| \leq \eta$, $\bar f=f(\bar{\bf x})$}.
\end{algorithmic}
\end{algorithm}

\begin{proposition}
  Algorithm~\ref{algo:3bis} terminates and $\bar{\bf x}$ is
  such that $h(\bar{\bf x}) \leq \frac{\eta}{\lambda^{\max}}$ and
  $|\bar f-f^*| \leq \eta$.
\end{proposition}

\begin{proof}
By contradiction, assume that the algorithm does not terminate.
Let $l_i$ be the value of $\lambda^{\min}$ returned by the $i$-th call
to \textbf{DualLagrangian}.
Sequence $l_i$ is monotone non increasing, moreover the domain of the
sequence is a subset of finite cardinality of interval $[0,
\lambda^{\tt init}]$ (its maximum cardinality depends on
$\lambda^{\tt init}$ and $\epsilon$). Indeed, the termination condition of function
\textbf{DualLagrangian} allows only for a finite number of divisions
of the interval $[0, \lambda^{\tt init}]$. 
Hence, sequence $l_i$ converges in a finite number of iterations to its
limit $l_\infty=\lim_{i \to \infty} l_i$ and there exists $\bar i \in
\mathbb{N}$ such that $(\forall i \geq \bar i)\; l_i=l_\infty$.
By iv) of Proposition~\ref{prop_alg}, $h^{\max}(l_\infty) \geq 0$ and,
since the algorithm does not terminate, $h^{\max}(l_\infty) \geq
\eta$.
At the $\bar i+1$-iteration, the algorithm calls
$\textbf{DualLagrangian}(X,l_\infty)$, which returns the value $\lambda^{\min}=l_\infty$.
Anyway, at the previous iteration $\bar i$, the elements $P_X(\lambda^{\min})$
at which function $h$ is positive had already been removed from $X$. This
implies that $\textbf{DualLagrangian}(X,l_\infty)$ cannot return the
strictly positive value $\lambda^{\min}=l_\infty$, leading to a contradiction.
Hence, the algorithm terminates and the stated bounds hold because
of the termination condition and by Proposition~\ref{prop_est}.
\end{proof}

\section{Lagrangian relaxation of the CDT problem}
\label{sec:dual}
In this section, we apply the general properties presented in
Section~\ref{sec_gen_prop} to the CDT problem~\eqref{eq:cdt}.
In fact, the CDT problem is a specific instance of~\eqref{eq:gen} in which $f({\bf x })=
 {\bf x }^\top {\bf Q} {\bf x}+{\bf q}^\top {\bf x}$, $g({\bf x })=
  {\bf x}^\top {\bf x}-1$, $h({\bf x })=
 {\bf x}^\top {\bf A}{\bf x}+{\bf a}^\top {\bf x}- a_0$.

\begin{figure}[!tbh]
\centering
\includegraphics[width=.9\columnwidth]{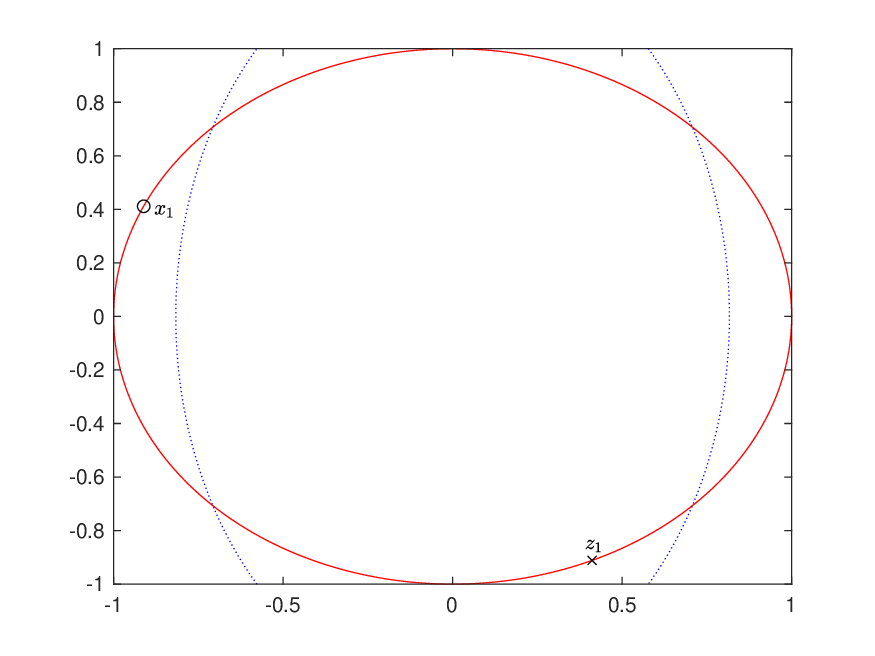}
\caption{Optimal solutions of the dual Lagrangian bound outside $H$ ($x_1$) and in $int(H)$ ($z_1$), denoted by $\circ$ and $\times$, respectively. The continuous red curve is the border of the unit ball, while the dotted blue curve is the border of the ellipsoid $H$.}
\label{fig:soldual}
\end{figure}

Note that the first two requirements of Assumption~\ref{assum_p} are
satisfied; in order to satisfy the third one we assume that
\begin{equation}
\label{eq:ella}
h_0=\min_{{\bf x}\ :\ {\bf x}^\top {\bf x}\leq 1} {\bf x}^\top {\bf A}{\bf x}+{\bf a}^\top {\bf x} -a_0<0,
\end{equation}
i.e., the feasible region of (\ref{eq:cdt}) has a nonempty interior. Note that the assumption can be checked in polynomial  time by the solution of a trust region problem.
%At the moment we do not assume that ${\bf A}$ is definite positive. Such assumption will be introduced later on.
As before, we denote by $X\subseteq \mathbb{R}^n$ a closed set such that $X
\supset H$, i.e., it contains the ellipsoid defined by the second constraint. 
For each $\lambda\geq 0$, the Lagrangian
relaxation~\eqref{eq:paramcdtg} takes on the form
\begin{equation}
\label{eq:paramcdt}
\begin{array}{ll}
p_X(\lambda)=\min_{{\bf x}\in X} & {\bf x }^\top ({\bf Q}+\lambda {\bf A}) {\bf x}+({\bf q}+\lambda {\bf a})^\top {\bf x} -\lambda a_0 \\ [6pt]
 &  {\bf x}^\top {\bf x}\leq 1.
\end{array}
\end{equation}
%As a consequence of Proposition~\ref{prop_bound}, for each $\lambda\geq 0$, the optimal value $p_X(\lambda)$ of this problem is a lower bound for problem (\ref{eq:cdt}). 
If $X=\mathbb{R}^n$, this is the standard Lagrangian relaxation of problem (\ref{eq:cdt}) and it can be solved efficiently since it is a trust region problem.
Following the notation of Section \ref{sec_gen_prop}, let
$$
P_X(\lambda)=\arg\min_{{\bf x}\in X\ :\ {\bf x}^\top {\bf x}\leq 1}  {\bf x }^\top ({\bf Q}+\lambda {\bf A}) {\bf x}+({\bf q}+\lambda {\bf a})^\top {\bf x} 
$$
be the set of optimal solutions of (\ref{eq:paramcdt}).
To apply Algorithm~\ref{algo:1} to the CDT problem with $X=\Real^n$,
we need to characterize the set of optimal solutions $P_{\mathbb{R}^n}(\lambda)$ of problem (\ref{eq:paramcdt}) with $X=\mathbb{R}^n$, which is a trust region problem.
The set of optimal solutions of a trust region problem has been derived, e.g., in \cite{Adachi17,More83,Sorensen82}. Here we briefly recall the different cases. 
For simplicity, let ${\bf S}_{\lambda}={\bf Q}+\lambda {\bf A}$ and ${\bf s}_{\lambda}={\bf q}+\lambda {\bf a}$. 
We distinguish the following cases:
\begin{description}
\item[Case 1] If ${\bf S}_{\lambda}\succ {\bf O}$ and$ \left\|-\frac{1}{2}{\bf S}_{\lambda}^{-1}{\bf s}_{\lambda}\right\|\leq 1$, then $-\frac{1}{2}{\bf S}_{\lambda}^{-1}{\bf s}_{\lambda}$ is the unique optimal solution of (\ref{eq:paramcdt});
\item[Case 2] Let ${\bf u}_j$ be the orthonormal eigenvectors of matrix  ${\bf S}_{\lambda}$, and let $\gamma_j$ be the corresponding eigenvalues. Let $\gamma_{\min}=\min_j \gamma_j$ and
$J_{\lambda}=\arg\min_j \gamma_j$. For each $\gamma$ such that $\forall j\ \gamma\neq \gamma_j$, let
$$
{\bf y}(\gamma)={\bf y}_1(\gamma)+{\bf y}_2(\gamma),
$$
where
$$
{\bf y}_1(\gamma)=-\sum_{j\not\in J_{\lambda}} \frac{{\bf s}_{\lambda}^\top {\bf u}_j}{\gamma_j-\gamma}{\bf u}_j,\ \ 
{\bf y}_2(\gamma)=-\sum_{j\in J_{\lambda}} \frac{{\bf s}_{\lambda}^\top {\bf u}_j}{\gamma_j-\gamma}{\bf u}_j.
$$
Then, we have the following subcases:
\begin{description}
\item[Case 2.1] It holds that  ${\bf s}_{\lambda}^\top {\bf u}_j\neq 0$ for some $j\in J_{\lambda}$. Then, 
there exists a unique $\gamma^*\in (-\gamma_{\min},+\infty)$ such that $\|{\bf y}(\gamma^*)\|=1$ and
${\bf y}(\gamma^*)$ is the unique optimal solution of (\ref{eq:paramcdt});
\item[Case 2.2] It holds that  ${\bf s}_{\lambda}^\top {\bf u}_j= 0$ for all $j\in J_{\lambda}$ but 
$\|{\bf y}_1(\gamma_{\min})\|\geq 1$. In this case there exists a unique $\gamma^*\in [-\gamma_{\min},+\infty)$ such that $\|{\bf y}_1(\gamma^*)\|=1$ and ${\bf y}_1(\gamma^*)$ is the unique optimal solution of (\ref{eq:paramcdt});
\item[Case 2.3] It holds that  ${\bf s}_{\lambda}^\top {\bf u}_j= 0$ for all $j\in J_{\lambda}$ and 
$\|{\bf y}_1(\gamma_{\min})\|< 1$. In this case
we have that $P_{\mathbb{R}^n}(\lambda)$ is not a singleton and is made up by the following points:
\begin{equation}
\label{eq:optimalset}
P_{\mathbb{R}^n}(\lambda)=\left\{{\bf y}_1(\gamma_{\min})+\sum_{j\in J_{\lambda}} \xi_j {\bf u}_j\ :\ \sum_{j\in J_{\lambda}} \xi_j^2=1-\|{\bf y}_1(\gamma_{\min})\|^2\right\}.
\end{equation}
Thus, we recognize two further subcases:
\begin{description}
\item[Case 2.3.1] $|J_{\lambda}|=1$, in which case $P_{\mathbb{R}^n}(\lambda)$ contains exactly two distinct points;
\item[Case 2.3.2] $|J_{\lambda}|\geq 2$, in which case the set $P_{\mathbb{R}^n}(\lambda)$ contains an infinite number of points and is a connected set.
\end{description}
\end{description}
\end{description}
Note that in Cases 2.3.1 and 2.3.2 we can compute the two values $h_{\mathbb{R}^n}^{\min}(\lambda),h_{\mathbb{R}^n}^{\max}(\lambda)$ by solving a trust region problem over the border of a $|J_{\lambda}|$-dimensional ball. More precisely, we need to solve the following problems:
\begin{equation}
\label{eq:hminmax}
\begin{array}{ll}
\min/\max_{\boldsymbol{\xi}} & {\bf w}(\boldsymbol{\xi})^\top {\bf A} {\bf w}(\boldsymbol{\xi}) + {\bf a}^\top {\bf w}(\boldsymbol{\xi})-a_0 \\ [6pt]
& \|{\bf w}(\boldsymbol{\xi})\|^2 = 1, %-\|{\bf y}_1(\gamma_{\min})\|^2,
\end{array}
\end{equation}
where ${\bf w}(\boldsymbol{\xi})={\bf y}_1(\gamma_{\min})+\sum_{j\in J_{\lambda}} \xi_j {\bf u}_j$. In these cases, where $P_{\mathbb{R}^n}(\lambda)$ is not a singleton, we also set 
\begin{equation}
\label{eq:defz1z2}
\begin{array}{ll}
{\bf z}_1(\lambda)={\bf w}(\boldsymbol{\xi}_1) & \boldsymbol{\xi}_1\in\arg\min_{\boldsymbol{\xi}\ :\ \|{\bf w}(\boldsymbol{\xi})\|=1}  {\bf w}(\boldsymbol{\xi})^\top {\bf A} {\bf w}(\boldsymbol{\xi}) + {\bf a}^\top {\bf w}(\boldsymbol{\xi})-a_0  \\ [8pt]
{\bf z}_2(\lambda)={\bf w}(\boldsymbol{\xi}_2) & \boldsymbol{\xi}_2\in\arg\max_{\boldsymbol{\xi}\ :\ \|{\bf w}(\boldsymbol{\xi})\|=1}  {\bf w}(\boldsymbol{\xi})^\top {\bf A} {\bf w}(\boldsymbol{\xi}) + {\bf a}^\top {\bf w}(\boldsymbol{\xi})-a_0,
\end{array}
\end{equation}
while
% is an optimal solution of the $\min$ problem in (\ref{eq:hminmax}), and 
%${\bf z}_2(\lambda)={\bf w}(\boldsymbol{\xi}_2)$, where ${\bf w}(\boldsymbol{\xi}_2)$ is an optimal solution of the $\max$ problem in (\ref{eq:hminmax}). %MA QUESTI NON POSSONO AVERE SOLUZIONI OTTIME MULTIPLE?????
in all other cases, when $P_{\mathbb{R}^n}(\lambda)=\{{\bf z}^\star(\lambda)\}$ is a singleton, we set
\begin{equation}
\label{eq:defz1z2singl}
{\bf z}_1(\lambda)={\bf z}_2(\lambda)={\bf z}^\star(\lambda).
\end{equation}
$\ $% are both set equal to the unique optimal solution of (\ref{eq:paramcdt}).
\newline\newline\noindent
%DIRE CHE DESCRIVIAMO ALGORITMO NEL CASO GENERALE MA CHE AL MOMENTO
%STIAMO IPOTIZZANDO $X=\mathbb{R}^n$.
%where $\hat{\lambda}$ is defined in~\eqref{eq:condlambda_g}.
The following statement is a direct consequence of
Proposition~\ref{prop:noexact_g}. %establishes when the bound returned by Algorithm~\ref{algo:1} is exact, i.e., equal to $p^\star$.
\begin{proposition}
\label{obs:jgeq2}
In the CDT Problem~\eqref{eq:cdt}, if
\begin{itemize}
\item $h_{\mathbb{R}^n}^{\min}(\lambda)=0$;
\item or if $h_{\mathbb{R}^n}^{\max}(\lambda)=0$;
\item or $h_{\mathbb{R}^n}^{\min}(\lambda)<0<h_{\mathbb{R}^n}^{\max}(\lambda)$ and $|J_{\lambda}|\geq 2$ (i.e., we are in Subcase 2.3.2);
\end{itemize}
then $p_{\mathbb{R}^n}(\lambda)=p^\star$.
\end{proposition}
\begin{proof}
Since $\{h_{\mathbb{R}^n}^{\min}(\lambda),
h_{\mathbb{R}^n}^{\max}(\lambda)\} \in Q_X(\lambda)$,
  in the first two cases $0 \in Q_X(\lambda)$ and the thesis is a
  consequence of
  Proposition~\ref{prop:noexact_g}.
%If  $h_{\mathbb{R}^n}^{\min}(\lambda)=0$, then there exists ${\bf y}\in P_{\mathbb{R}^n}(\lambda)\cap \partial {\cal E}$ which is also feasible for (\ref{eq:cdt}) and such that
%$$
%p_{\mathbb{R}^n}(\lambda)={\bf y }^\top {\bf Q} {\bf y}+{\bf q}^\top {\bf y},
%$$
%so that ${\bf y}$ is, in fact, optimal for problem (\ref{eq:cdt}) and $p_{\mathbb{R}^n}(\lambda)$ is the optimal value of this problem.
%\newline\newline\noindent
%The case $h_{\max}^{\lambda}(\mathbb{R}^n)=0$ is completely similar.  
%\newline\newline\noindent
If $h_{\mathbb{R}^n}^{\min}(\lambda)<0< h_{\mathbb{R}^n}^{\max}(\lambda)$ and $|J_{\lambda}|\geq2$, we observed that  $P_{\mathbb{R}^n}(\lambda)$ is a connected set. Then, there exists ${\bf x}^\star \in P_{\mathbb{R}^n}(\lambda)$ such that
${\bf x}^\star\in \partial H$. More precisely, ${\bf x}^\star$ is a point along the curve in $P_{\mathbb{R}^n}(\lambda)$ connecting points ${\bf z}_1(\lambda)$ and ${\bf z}_2(\lambda)$, defined in (\ref{eq:defz1z2}). Thus, the lower bound $p_{\mathbb{R}^n}(\lambda)$ is equal to the optimal value of problem (\ref{eq:cdt}).
\end{proof}
Note that the first two conditions of Proposition~\ref{obs:jgeq2} imply exactness of the bound also for generic regions $X\supset H$, while the last condition is specific to the case $X=\mathbb{R}^n$.
%If we remove the exit condition $\lambda^X_{\max}-\lambda^X_{\min}>\varepsilon$, then either the algorithm stops after a finite number of iterations, as soon as the $exit$ flag is set equal to ${\tt true}$, or it runs for an infinite number of iterations
%and $\lambda^X_{\min},\lambda^X_{\max}\rightarrow \bar{\lambda}^X$. For such limit value we have that
%$$
%H_{\min}^{\bar{\lambda}^X}=\lim H_{\max}^{\lambda_{\max}^X},\ \ \ H_{\max}^{\bar{\lambda}^X}=\lim H_{\min}^{\lambda_{\min}^X},
%$$
%If, after setting $\varepsilon=0$, Algorithm \ref{algo:1} does not stop in a finite number of iterations and $\lambda^{\mathbb{R}^n}_{\min},\lambda^{\mathbb{R}^n}_{\max}\rightarrow \bar{\lambda}^{\mathbb{R}^n}$,
%then, as for the general case, we have
%$0\in \left[h_{\min}^{\bar{\lambda}^{\mathbb{R}^n}}(\mathbb{R}^n),h_{\max}^{\bar{\lambda}^{\mathbb{R}^n}}(\mathbb{R}^n)\right]$, so that $\bar{\lambda}^{\mathbb{R}^n}$ is the maximizer of $p_{\mathbb{R}^n}$ and $p_{\mathbb{R}^n}(\bar{\lambda}^{\mathbb{R}^n})$ is the dual Lagrangian bound. In particular, if $h_{\min}^{\bar{\lambda}^{\mathbb{R}^n}}(\mathbb{R}^n)=0$, or $h_{\max}^{\bar{\lambda}^{\mathbb{R}^n}}(\mathbb{R}^n)=0$, or if at $\bar{\lambda}^{\mathbb{R}^n}$ it holds that $|J_{\bar{\lambda}^{\mathbb{R}^n}}|\geq 2$ (see Observation \ref{obs:jgeq2}), then we can conclude that $p_{\mathbb{R}^n}(\bar{\lambda}^{\mathbb{R}^n})=p^\star$, i.e., the dual Lagrangian bound is exact. 
%We have thus proved the following result, related 
%Thus, the only case when the dual Lagrangian bound is not exact is 
The following result is related
to the necessary and sufficient condition under which the dual Lagrangian bound is not exact discussed in~\cite{Ai09}.
\begin{proposition}
  \label{prop:noexact}
  In the CDT problem~\eqref{eq:cdt},
  $p_{\mathbb{R}^n}(\lambda_{\mathbb{R}^n}) \neq p^*$ if and only if
$P_{\mathbb{R}^n}(\lambda_{\mathbb{R}^n})$ contains exactly two points, i.e., $|J_{\lambda_{\mathbb{R}^n}}|=1$ (Subcase 2.3.1), and $0\in \left(h_{\mathbb{R}^n}^{\min}(\lambda_{\mathbb{R}^n}), h_{\mathbb{R}^n}^{\max}(\lambda_{\mathbb{R}^n})\right)$.
\end{proposition}
\begin{proof}
  It is a consequence of Proposition~\ref{obs:jgeq2} and the fact that for $|J_{\lambda_{\mathbb{R}^n}}|=1$   it holds that $Q_{\mathbb{R}^n}(\lambda_{\mathbb{R}^n})=\{h_{\mathbb{R}^n}^{\min}(\lambda_{\mathbb{R}^n}), h_{\mathbb{R}^n}^{\max}(\lambda_{\mathbb{R}^n})\}\not\ni 0$.
\end{proof}
Now, we introduce an example where
$p_{\mathbb{R}^n}(\lambda_{\mathbb{R}^n}) \neq p^*$, that is the dual Lagrangian bound is not exact, which will also be helpful in the following sections.
%We illustrateThis is also strictly related  to the necessary and sufficient condition under which the dual Lagrangian bound is exact discussed in \cite{Ai09}.
\begin{example}
\label{ex:1}
Let us consider the following example taken from~\cite{Burer13}:
$$
{\bf Q}=\left(\begin{array}{cc}
-4 & 1 \\ 
1 & -2
\end{array}
\right),\ \ {\bf A}=\left(\begin{array}{cc}
3 & 0 \\
0 & 1
\end{array}
\right),\ \  {\bf q}=\left(1\ 1\right)\ \ {\bf a}=\left(0\ 0\right),\ \ a_0=2.
$$
Such instance has optimal value $-4$ attained at points $\left(\frac{\sqrt{2}}{2},-\frac{\sqrt{2}}{2}\right)$ and $\left(-\frac{\sqrt{2}}{2}, \frac{\sqrt{2}}{2}\right)$. The maximizer of $p_{\mathbb{R}^2}(\lambda)$ is
$\lambda_{\mathbb{R}^2}=1$ for which we have:
$$
h_{\mathbb{R}^2}^{\min} \approx -0.66  <0< 0.66 \approx h_{\mathbb{R}^2}^{\max},
$$
and, moreover, $|J_{\lambda_{\mathbb{R}^2}}|=1$, so that we have exactly two optimal solutions of (\ref{eq:paramcdt}), one violating the second constraint, namely ${\bf z}_2(\lambda_{\mathbb{R}^2})=(-0.911, 0.4114)$, 
point $x_1$ in Figure \ref{fig:soldual}, displayed as $\circ$, the other  in $int(H)$, point $z_1$ in Figure \ref{fig:soldual}, displayed as $\times$. The lower bound is $p_{\mathbb{R}^2}(1)=-4.25$, which is not exact. 
\end{example}
In the next sections we will try to improve the dual Lagrangian bound (or the equivalent SDP bound) by adding linear cuts, i.e., by introducing regions $X$ defined by one or two linear cuts.
\section{Bound improvement}
\label{sec:boundimp}
\begin{figure}[!tbh]
\centering
\includegraphics[width=.9\columnwidth]{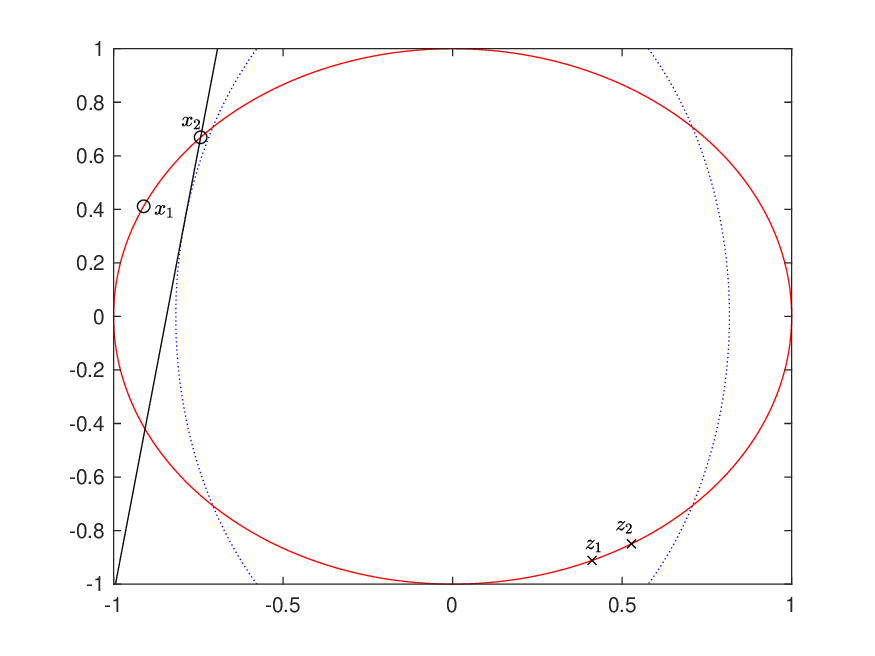}
\caption{First linear cut and the two optimal solutions lying outside $H$ ($x_2$) and in $int(H)$ ($z_2$), denoted by $\circ$ and $\times$, respectively.}
\label{fig:solfirstlincut}
\end{figure}

%In this section we introduce the additional assumption that ${\bf A}\succ {\bf O}$, i.e., the second constraint is also a strictly convex one. 
%We will denote by ${\cal E}$ the ellipsoid defined by the second constraint, and by $\partial {\cal E}$ its border.
We assume that the dual Lagrangian relaxation is not exact, i.e., as previously stated in Proposition \ref{prop:noexact}
%there exists $\bar{\lambda}^{\mathbb{R}^n}$ such that
$$
0\in \left(h_{\mathbb{R}^n}^{\min}(\lambda_{\mathbb{R}^n}), h_{\mathbb{R}^n}^{\max}(\lambda_{\mathbb{R}^n})\right),\ \ \ |J_{\lambda_{\mathbb{R}^n}}|=1.
$$
Recall that, by Proposition \ref{prop:noexact}, in this case, there exists a single point ${\bf z}_1(\lambda_{\mathbb{R}^n})\in P_{\mathbb{R}^n}(\lambda_{\mathbb{R}^n})\cap H$ (actually ${\bf z}_1(\lambda_{\mathbb{R}^n})\in int(H)$), and a single point  ${\bf z}_2(\lambda_{\mathbb{R}^n})\in P_{\mathbb{R}^n}(\lambda_{\mathbb{R}^n})\setminus H$.
Now we show that the dual Lagrangian bound can be strictly improved through the addition of a linear cut.
We first observe that  the optimal value of  problem (\ref{eq:cdt}) does not change if we 
add constraints which are implied by the second one.
\newline\newline\noindent
In the following proposition, we define a  projection $\Pi_{{\bf A},{\bf a}}:
\Real^n \setminus H \to \partial H$, that
maps ${\bf x}\not\in H$ to the element of $\partial H$ located
on the segment that joins ${\bf x}$ to the center of the ellipsoid $H$ (given by $\boldsymbol{\alpha}=-\frac{1}{2}
{\bf A}^{-1} {\bf a}$).
\begin{proposition}
  For ${\bf x} \not \in H$,
  set $\Pi_{{\bf A}, {\bf a}}({\bf x})= \sqrt{\frac{-h(\boldsymbol{\alpha})}{h({\bf x})-h(\boldsymbol{\alpha})}} ({\bf x}-\boldsymbol{\alpha}) + \boldsymbol{\alpha}$, where $\boldsymbol{\alpha}=-\frac{1}{2}
{\bf A}^{-1} {\bf a}$ is the center of the ellipsoid.
Then $h(\Pi_{{\bf A},{\bf a}}({\bf x}))=0$.
\end{proposition}
\begin{proof}
Note that 
$(\forall \beta \in \Real) \; h(\beta({\bf x}-\boldsymbol{\alpha})+\boldsymbol{\alpha}) -
h(\boldsymbol{\alpha})=\beta^2 (h({\bf x})-h(\boldsymbol{\alpha}))$ (it is a consequence of the fact
that function $h$ is quadratic and it can be verified by direct substitution).
Then $h(\Pi_{{\bf A},{\bf a}}({\bf x}))=h\left( \sqrt{\frac{-h(\boldsymbol{\alpha})}{h({\bf x})-h(\boldsymbol{\alpha})}} ({\bf x}-\boldsymbol{\alpha}) +
\boldsymbol{\alpha}\right)=\frac{-h(\boldsymbol{\alpha})}{h({\bf x})-h(\boldsymbol{\alpha})} (h({\bf x})-h(\boldsymbol{\alpha})) + h(\boldsymbol{\alpha})=0$.
\end{proof}
Given any $\bar{{\bf x}} \in \Real^n$, it holds, by convexity, that
$$
{\bf x}^\top {\bf A}{\bf x}+{\bf a}^\top {\bf x}\geq   \bar{{\bf x}}^\top {\bf A}\bar{{\bf x}}+{\bf a}^\top \bar{{\bf x}}+(2 {\bf A}\bar{{\bf x}}+{\bf a})^\top({\bf x}-\bar{{\bf x}}).
$$
Thus, the following linear constraint is implied by the second constraint in (\ref{eq:cdt}):
%one can add the following linear constraint (supporting hyperplane at $\bar{{\bf x}}$):
\begin{equation}
\label{eq:linearcut}
(2 {\bf A}\bar{{\bf x}}+{\bf a})^\top {\bf x}-\bar{{\bf x}}^\top {\bf A} \bar{{\bf x}} \leq a_0,
\end{equation}
and, consequently, it can be added to problem (\ref{eq:cdt}) without
modifying its feasible region.
In particular, if $\bar {\bf x} \in \partial H$, being $\bar {\bf x}^T {\bf A} \bar {\bf x} + {\bf a}^T \bar {\bf x}=a_0$,
the linear constraint is
$$
(2{\bf A} \bar{{\bf x}}+{\bf a})^\top ({\bf x}-\bar{{\bf x}})\leq 0.
$$
%
%Thus, the following linear constraint is implied by the second constraint in (\ref{eq:cdt}):
%\begin{equation}
%\label{eq:linearcut}   (2 {\bf
%  A}\bar{{\bf x}}+{\bf a})^\top({\bf x}-\bar{{\bf x}}) \leq 0,
%\end{equation}
%and, consequently, it can be added to problem (\ref{eq:cdt}) without
%modifying its feasible region.
Due to the redundancy of the lienar constraint for problem (\ref{eq:cdt}), we can define, for a given $\bar {\bf x} \in \partial H$,  the new Lagrangian problem
\begin{equation}
\label{eq:paramcdt1}
\begin{array}{ll}
p_X(\lambda)=\min_{{\bf x}} & {\bf x }^\top ({\bf Q}+\lambda {\bf A}) {\bf x}+({\bf q}+\lambda {\bf a})^\top {\bf x} -\lambda a_0 \\ [6pt]
 &  {\bf x}^\top {\bf x}\leq 1 \\ [6pt]
& (2{\bf A} \bar{{\bf x}}+{\bf a})^\top ({\bf x}-\bar{{\bf x}})\leq 0.
\end{array}
\end{equation}
where
%Note that this is exactly problem (\ref{eq:paramcdt}) where
\begin{equation}
\label{eq:hypersp}
X=\Omega_{\bar{{\bf x}}}=\{{\bf x}\ :\ (2 {\bf A}\bar{{\bf x}}+{\bf a})^\top({\bf x}-\bar{{\bf x}})\leq 0\}\supset H.
\end{equation}
%Theorem \ref{lem:1} still holds true. We recall definition (\ref{eq:defdiffh}) of the set $H_{\Omega_{\bar{\bf x}}}^{\lambda}$ and of the related values $h_{\Omega_{\bar{\bf x}}}^{\min}(\lambda)$ and $h_{\Omega_{\bar{\bf x}}}^{\max}(\lambda)$.
If we set  $\bar{\bf  x}=\Pi_{{\bf A},{\bf a}}
({\bf z}_2(\lambda_{\mathbb{R}^n})) $, i.e., $\bar{\bf  x}$ is the projection over $\partial H$ of 
the single point in $P_{\mathbb{R}^n}(\lambda_{\mathbb{R}^n})\setminus H$, then $\mathbb{R}^n \supset X \supset H$ and, moreover, $X\cap(P_{\mathbb{R}^n}(\lambda_{\mathbb{R}^n})\setminus H)=\emptyset$, so that, by Proposition
\ref{prop:boundimprove1}, $\bar{p}_X> \bar{p}_{\mathbb{R}^n}$.
Then, if we run again Algorithm \ref{algo:1} with input $X=\Omega_{\bar{\bf x}}$ defined in (\ref{eq:hypersp}) and $\lambda_{\tt init}=\lambda_{\mathbb{R}^n}$ (or $\lambda_{\tt init}={\lambda}^{\max}_{\mathbb{R}^n}$), we are able 
to improve strictly the dual Lagrangian bound.
Note that
 problem (\ref{eq:paramcdt1}), needed to compute function $p_{\Omega_{\bar{{\bf x}}}}$, can be solved in polynomial time according to the results proved in \cite{Burer13,Sturm03}.
But we also discuss an alternative way to solve problem (\ref{eq:paramcdt1}), based on the solution of a trust region problem.
%implement Algorithm \ref{algo:1}, namely Algorithm \ref{algo:2}, which will turn out to be useful in what follows and requires to solve a trust region problem at each iteration. 
For $\lambda=\lambda_{\mathbb{R}^n}$, after the addition of the linear cut, a unique optimal solution exists, lying in 
$int(H)$  and, consequently, in $int(\Omega_{\bar{{\bf x}}})$, since also the linear constraint in (\ref{eq:paramcdt1}) is not active at it, being $H$ a subset of the region defined by the linear cut. By continuity, for $\lambda$ values smaller than but close to $\lambda_{\mathbb{R}^n}$, 
the unique optimal solution of  (\ref{eq:paramcdt1}) also lies in $int(H)$, i.e., $P_{\Omega_{\bar{\bf x}}}(\lambda)=\{{\bf z}_1(\lambda)\}$ with ${\bf z}_1(\lambda)\in int(H)$. Thus, such optimal solution must be a local and nonglobal optimal solution
of the trust region problem (\ref{eq:paramcdt}) with $X=\mathbb{R}^n$. Indeed, the globally optimal solution of this trust region problem always violates the second constraint in (\ref{eq:cdt}) for all $\lambda<\lambda_{\mathbb{R}^n}$.
Now, for a generic $\lambda\in [0,\lambda_{\mathbb{R}^n})$, we first check
%At this point we can run Algorithm \ref{algo:2}. For each tested $\lambda$ value in the while loop we can first check 
whether a local and nonglobal optimal solution of problem (\ref{eq:paramcdt}) with $X=\mathbb{R}^n$ exists, by exploiting the necessary and sufficient condition stated in \cite{WangXia20}. Also recall that, if it exists, the local and nonglobal minimizer is unique. If it does not exist, then 
the linear constraint must be active at all optimal solutions of problem (\ref{eq:paramcdt1}). In this case we set $f_1=+\infty$. Otherwise, if it exists, we denote it by ${\bf z}_1(\lambda)$. If  ${\bf z}_1(\lambda)\not \in \Omega_{\bar{{\bf x}}}$, then we set again $f_1=+\infty$, otherwise we denote by $f_1$ the value of the objective function of (\ref{eq:paramcdt1}) evaluated at ${\bf z}_1(\lambda)$.
Then, we consider the best feasible solutions  of problem (\ref{eq:paramcdt1}) for which the linear constraint is imposed to be active. The resulting problem is converted into a trust region problem, after
the change of variable ${\bf x}=\bar{{\bf x}}+{\bf V}{\bf z}$, where ${\bf V}\in \mathbb{R}^{n\times (n-1)}$ is a matrix whose columns form a basis for the null space of vector $2{\bf A}\bar{{\bf x}}+{\bf a}$. 
The resulting (trust region) problem is:
\begin{equation}
\label{eq:trustnull}
\begin{array}{ll}
\min_{{\bf w} \in \mathbb{R}^{n-1}} & {\bf w}^\top {\bf V}^\top ({\bf Q}+\lambda {\bf A}) {\bf V} {\bf w}+\left[2 \bar{{\bf x}}^\top  ({\bf Q}+\lambda {\bf A}) {\bf V} + ({\bf q}+\lambda {\bf a})^\top\right] {\bf w} +\ell(\bar{{\bf x}},\lambda)
 \\ [6pt]
& \| \bar{{\bf x}}+{\bf V} {\bf w}\|^2 \leq 1,
\end{array}
\end{equation}
where $\ell(\bar{{\bf x}},\lambda)=
\bar{{\bf x}}^\top  ({\bf Q}+\lambda {\bf A})\bar{{\bf x}}+
({\bf q}+\lambda {\bf a})^\top \bar{{\bf x}} - \lambda a_0$ is constant with respect to the vector of variables ${\bf w}$.
Let $W^\star(\lambda)$ be the set of optimal solutions of problem (\ref{eq:trustnull}) and 
$$
P_1^\star(\lambda)=\{ \bar{{\bf x}}+{\bf V} {\bf w}^\star\ : \ {\bf w}^\star\in W^\star(\lambda)\}. 
$$
Note that the set $W^\star(\lambda)$ can be computed through the procedure presented in Section \ref{sec:dual} with the different cases (namely, Cases 1, 2.1, 2.2, 2.3.1, 2.3.2) after rewriting it as a classical trust region problem.
Moreover, let $f_2<+\infty$ be the optimal value of problem (\ref{eq:trustnull}). Now, after comparing $f_1$ and $f_2$, we are able to define the set $P_{\Omega_{\bar{\bf x}}}(\lambda)$ of optimal solutions for problem (\ref{eq:paramcdt1}).
%Note that Then, we set ${\bf z}_2(\lambda)=\bar{{\bf x}}+{\bf V} {\bf w}^\star$ and denote by $f_2$ its objective function value.
More precisely, if $f_2>f_1$, then $P_{\Omega_{\bar{\bf x}}}(\lambda)=\{{\bf z}_1(\lambda)\}$, i.e.,  ${\bf z}_1(\lambda)$ is the unique optimal solution of problem (\ref{eq:paramcdt1}). 
In this case 
$$
h_{\Omega_{\bar{\bf x}}}^{\min}(\lambda)=h_{\Omega_{\bar{\bf x}}}^{\max}(\lambda)={\bf z}_1(\lambda)^\top {\bf A} {\bf z}_1(\lambda)+{\bf a}^\top {\bf z}_1(\lambda)-a_0.
$$
%Then, if ${\bf z}_1(\lambda)$ lies in the interior of ${\cal E}$, then we set  $\lambda_{\max}=\lambda$; if ${\bf z}_1(\lambda)\in \partial {\cal E}$, then ${\bf z}_1(\lambda)$ is an optimal solution of the original problem (\ref{eq:cdt}) and the bound is exact; finally, if ${\bf z}_1(\lambda)\not\in {\cal E}$, we set $\lambda_{\min}=\lambda$.
Instead, if $f_2<f_1$, which always holds, e.g., if $f_1=+\infty$, then $P_{\Omega_{\bar{\bf x}}}(\lambda)=P_1^\star(\lambda)$. Since all points in $P_1^\star(\lambda)$ lie over a supporting hyperplane of $H$, we must have that
$$
h_{\Omega_{\bar{\bf x}}}^{\min}(\lambda)=\min_{ {\bf x}\in P_1^\star(\lambda)} {\bf x}^\top {\bf A} {\bf x}+{\bf a}^\top {\bf x}-a_0  \geq 0,
$$
and equality holds only if $\bar{{\bf x}}\in P_1^\star(\lambda)$. In the latter case, the bound is exact, otherwise Algorithm \ref{algo:1} sets $\lambda^{\min}=\lambda$. 
%Otherwise, we set $\lambda_{\min}=\lambda$. 
Finally,
if $f_1=f_2$, then $P_{\Omega_{\bar{\bf x}}}(\lambda)=P_1^\star(\lambda) \cup \{{\bf z}_1(\lambda)\}$ and in this case $0\in [h_{\Omega_{\bar{\bf x}}}^{\min}(\lambda), h_{\Omega_{\bar{\bf x}}}^{\max}(\lambda)]$ and the algorithms exits the loop.
%Then: if ${\bf z}_1(\lambda)\in \partial {\cal E}$, then ${\bf z}_1(\lambda)$ is an optimal solution of the original problem (\ref{eq:cdt}) and the bound is exact; if ${\bf z}_1(\lambda)\not\in {\cal E}$, we set $\lambda_{\min}=\lambda$;  if $\bar{{\bf x}}\in F_1^\star(\lambda)$, then the bound is exact; otherwise (${\bf z}_1(\lambda)$ lies in the interior of ${\cal E}$ and $\bar{{\bf x}}\not \in F_1^\star(\lambda)$,
%then $h_{\min}^{\lambda}(X)<0<h_{\max}^{\lambda}(X)$ and the algorithm stops.
The following result is a straightforward consequence Proposition
\ref{prop:boundimprove1}.
\begin{proposition}
\label{prop:improve}
Algorithm \ref{algo:1}  with $\varepsilon=0$ will stop after a finite number of iterations or will converge to some $\lambda_{\Omega_{\bar{\bf x}}}<\lambda_{\mathbb{R}^n}$ 
with a new lower bound $\bar{p}_{\Omega_{\bar{\bf x}}}>\bar{p}_{\mathbb{R}^n}$.
\end{proposition}
\begin{proof}
Strict inequalities hold in view of
 Proposition
\ref{prop:boundimprove1} with $X=\mathbb{R}^n$ and $Y=\Omega_{\bar{\bf x}}$, since, as already observed,
 $\Omega_{\bar{\bf x}}\cap(P_{\mathbb{R}^n}(\lambda_{\mathbb{R}^n})\setminus H)=\emptyset$.
\end{proof}
If the final bound is not exact, i.e., $\bar{p}_{\Omega_{\bar{\bf x}}}=p_{\Omega_{\bar{\bf x}}}(\lambda_{\Omega_{\bar{\bf x}}})< p^\star$, at $\lambda_{\Omega_{\bar{\bf x}}}$ we have
$f_1=f_2$ and $P_{\Omega_{\bar{\bf x}}}(\lambda_{\Omega_{\bar{\bf x}}})$ contains multiple optimal solutions, in particular, one in $int(H)$ and the other(s) outside
$H$. %	QUI CAMBIARE E RICONTROLLARE TUTTO ALGORITMO PERCHE' RESTITUSICA L'INSIEME DI TUTTE LE SOLUZIONI. namely vector ${\bf v}$ returned by Algorithm \ref{algo:2}.
We illustrate all this on Example \ref{ex:1}. 
\begin{example}
The optimal solution of (\ref{eq:paramcdt}) with $X=\mathbb{R}^n$ for
$\lambda_{\mathbb{R}^n}=1$ which violates the second constraint is
${\bf z}_2(\lambda_{\mathbb{R}^n})=(-0.911, 0.4114)$. The lower bound
is $p_{\mathbb{R}^n}(1)=-4.25$. After the addition of the linear
inequality (\ref{eq:linearcut}) obtained with $\bar{{\bf x}}=
\Pi_{{\bf A},{\bf a}}({\bf z}_2(\lambda_{\mathbb{R}^n}))$, equal to the projection of  ${\bf z}_2(\lambda_{\mathbb{R}^n})$ over the boundary of the second constraint, we can run again
Algorithm~\ref{algo:1} with $X=\Omega_{\bar{{\bf x}}}$ and we get to $\lambda_{\Omega_{\bar{{\bf x}}}}\approx 0.726$ and $p_{\Omega_{\bar{{\bf x}}}}(\lambda_{\Omega_{\bar{{\bf x}}}})\approx -4.097$, which improves the previous lower bound. 
In Figure \ref{fig:solfirstlincut} we show the linear cut and the two new optimal solutions outside $H$ and in 
$int(H)$ ($x_2$ and $z_2$, respectively) obtained at $\lambda_{\Omega_{\bar{{\bf x}}}}$. In the same figure we also display the previous pair of optimal solutions in order to show the progress of the algorithm.
\end{example}

%\begin{algorithm}
%\caption{Bound improvement through the additional linear inequality (\ref{eq:linearcut}).}
%\label{algo:2}
%\begin{algorithmic}
%\STATE{Set $Lb=-\infty$}
%\STATE{Set $\lambda_{\min}=0$, $\lambda_{\max}=\bar{\lambda}$}
%\STATE{Let $V$ be a matrix whose columns form a basis for the null space of vector $2{\bf A}\bar{{\bf x}}+{\bf a}$}
%\WHILE{$\lambda_{\max}-\lambda_{\min}>\varepsilon$}
%\STATE{Set $\lambda=(\lambda_{\max}+\lambda_{\min})/2$}
%\IF{Problem  (\ref{eq:paramcdt}) has no local and nonglobal minimizer}
%\STATE{Set $f_1=+\infty$} 
%\ELSE\STATE{Let ${\bf z}_1(\lambda)$ be the local and nonglobal minimizer for problem (\ref{eq:paramcdt})}
%\STATE{Set $f_1={\bf z}_1(\lambda)^\top ({\bf Q}+\lambda {\bf A}){\bf z}_1(\lambda)+({\bf q}+\lambda{\bf a})^\top {\bf z}_1(\lambda)$}
%\ENDIF
%\STATE{Let ${\bf z}_2(\lambda)=\bar{{\bf x}}+{\bf V} {\bf w}^\star$, where ${\bf w}^\star$ is a global minimzer of problem (\ref{eq:trustnull})}
%\STATE{Set $f_2={\bf z}_2(\lambda)^\top ({\bf Q}+\lambda {\bf A}){\bf z}_2(\lambda)+({\bf q}+\lambda{\bf a})^\top {\bf z}_2(\lambda) $}
%\IF{$f_2<f_1$ \OR ${\bf z}_1(\lambda)\not\in {\cal E}$}\STATE{Set $\lambda_{\min}=\lambda$, ${\bf v}={\bf z}_2(\lambda)$}
%\ELSE\STATE{Set $\lambda_{\max}=\lambda$}
%\ENDIF
%\IF{$\min\{f_1,f_2\}>Lb$}\STATE{Set $Lb=\min\{f_1,f_2\}$}\ENDIF
%\ENDWHILE
%\RETURN ${\bf v}, Lb, \lambda$
%\end{algorithmic}
%\end{algorithm}
%$\ $\newline\newline\noindent
\section{Improving the bound by local adjustments of the linear cut}
%\label{sec:further}
%We propose two different ways to further improve the bound. The first one is still based on the addition of a single linear cut but tries to adjust locally the point $\bar{{\bf x}}$ employed in (\ref{eq:paramcdt1}). The second one is based on the addition of a further linear cut.
%\subsection{Local adjustment of point $\bar{{\bf x}}$}
\label{sec:optimal}
In the previous section we proposed to set $\bar{{\bf x}}$ equal to the projection over $\partial H$ of  ${\bf z}_2(\lambda_{\mathbb{R}^n})$, the optimal solution of problem (\ref{eq:paramcdt})  with $X=\mathbb{R}^n$ lying outside $H$. However, this point can be improved by some local adjustment. To this end, we should search for some perturbation direction ${\bf d}$ such that
$$
[2 {\bf A}(\bar{{\bf x}}+{\bf d})+{\bf a}]^\top {\bf v}  -(\bar{{\bf x}}+{\bf d}) ^\top {\bf A} (\bar{{\bf x}}+{\bf d}) > a_0,
$$
for all ${\bf v}\in P_1^\star(\lambda_{\Omega_{\bar{{\bf x}}}})$.
% Indeed, if such a direction ${\bf d}$ exists, the linear cut (\ref{eq:linearcut}) is violated by all optimal solutions at which the current linear cut is active and, as we will see, by perturbing $\bar{\bf x}$ along this direction, we are able to further reduce $\lambda$ and, thus, to improve the bound.
Taking into account that the linear cut in (\ref{eq:paramcdt1}) is active at $\bar{{\bf x}}$, the above inequality is equivalent to
\begin{equation}
\label{eq:dirimp}
-{\bf d} ^\top {\bf A} {\bf d}+ 2 {\bf d}^\top {\bf A} ({\bf v}- \bar{{\bf x}})>0\ \ \ \forall {\bf v}\in  P_1^\star(\lambda_{\Omega_{\bar{{\bf x}}}}).
\end{equation}
If such direction ${\bf d}$ exists, then we are able to improve the linear cut by replacing $\bar{{\bf x}}$ with 
\begin{equation}
\label{Eq:defU}
\tilde{{\bf x}}=\Pi_{{\bf A},{\bf a}}(\bar{{\bf x}}+\eta {\bf d})\in \partial H,
\end{equation}
for some $\eta>0$ and small enough. 
%In fact, point $\tilde{{\bf x}}$ is the optimal solution of problem (\ref{eq:support}) where $\bar{{\bf x}}$ is replaced by $\bar{{\bf x}}+\eta {\bf d}$ (which can be obtained in closed form).
By (\ref{eq:dirimp}), for any positive step $\eta$ along direction ${\bf d}$ we have $P_1^\star(\lambda_{\Omega_{\bar{{\bf x}}}})\cap \Omega_{\tilde{{\bf x}}}=\emptyset$ and, by continuity, that holds true also in a small neighborhood
of $P_1^\star(\lambda_{\Omega_{\bar{{\bf x}}}})$. However, if we take too large a step along direction ${\bf d}$, then new optimal solutions of problem (\ref{eq:paramcdt1}) with $\bar{{\bf x}}$ replaced by $\tilde{{\bf x}}$, sufficiently far 
from $P_1^\star(\lambda_{\Omega_{\bar{{\bf x}}}})$, may appear. But if we take a small enough step along direction ${\bf d}$, then no new optimal solution will appear and the only optimal solution of problem (\ref{eq:paramcdt1}) with $\bar{{\bf x}}$ replaced by $\tilde{{\bf x}}$
will be point ${\bf z}_1(\lambda_{\Omega_{\bar{{\bf x}}}})\in int(H)$. Thus, 
as a consequence of Proposition \ref{prop:boundimprove},
Algorithm \ref{algo:1} with input $X=\Omega_{\tilde{{\bf x}}}$ will be able to further reduce the value $\lambda$ and improve (increase) the lower bound.
%choice is not guaranteed to be the one which leads to the best possible bound.
Therefore, the question now is how to find a direction ${\bf d}$ fulfilling (\ref{eq:dirimp}) or to establish it does not exist. We discuss different cases depending on the cardinality of $P_1^\star(\lambda_{\Omega_{\bar{{\bf x}}}})$ (see the cases discussed in Section \ref{sec:dual} for the trust region problem).
\subsection{$|P_1^\star(\lambda_{\Omega_{\bar{{\bf x}}}})|=1$}
\label{sec:x1=1}
In this case, let ${\bf v}$ be the unique point in $P_1^\star(\lambda_{\Omega_{\bar{{\bf x}}}})$, then we need to solve the following convex optimization problem
$$
\max_{{\bf d}\in \mathbb{R}^n} -{\bf d} ^\top {\bf A} {\bf d}+ 2 {\bf d}^\top {\bf A} ({\bf v}- \bar{{\bf x}}),
$$
whose optimal solution is ${\bf d}={\bf v}- \bar{{\bf x}}$ and its optimal value is $({\bf v}- \bar{{\bf x}})^\top {\bf A} ({\bf v}- \bar{{\bf x}})>0$.
Therefore, if $|P_1^\star(\lambda_{\Omega_{\bar{{\bf x}}}})|=1$ we are always able to locally adjust the current point $\bar{{\bf x}}$ in such a way that the bound can be improved.
\subsection{$|P_1^\star(\lambda_{\Omega_{\bar{{\bf x}}}})|=2$}
\label{sec:x1=2}
In this case, let ${\bf v}_1$ and ${\bf v}_2$ be the two optimal points in $P_1^\star(\lambda_{\Omega_{\bar{{\bf x}}}})$. Then, we need to solve the following optimization problem
\begin{equation}
\label{eq:optdir2}
\max_{{\bf d}\in \mathbb{R}^n} \min\{-{\bf d} ^\top {\bf A} {\bf d}+ 2 {\bf d}^\top {\bf A} ({\bf v}_1- \bar{{\bf x}}),-{\bf d} ^\top {\bf A} {\bf d}+ 2 {\bf d}^\top {\bf A} ({\bf v}_2- \bar{{\bf x}})\},
\end{equation}
or, equivalently
$$
\begin{array}{ll}
\max & v \\ [6pt]
 & v\leq  -{\bf d} ^\top {\bf A} {\bf d}+ 2 {\bf d}^\top {\bf A} ({\bf v}_1- \bar{{\bf x}}) \\ [6pt]
 & v\leq  -{\bf d} ^\top {\bf A} {\bf d}+ 2 {\bf d}^\top {\bf A} ({\bf v}_2- \bar{{\bf x}}).
\end{array}
$$
This is a convex optimization problem, whose solution can be obtained in closed form. Indeed, by imposing the KKT conditions, it can be seen that
the optimal solution has the following form
\begin{equation}
\label{eq:optform}
{\bf d}=\beta ({\bf v}_1- \bar{{\bf x}}) +(1-\beta)({\bf v}_2- \bar{{\bf x}}), \ \ \ \beta\in [0,1].
\end{equation}
Now, let 
$$
\begin{array}{l}
a=({\bf v}_1- \bar{{\bf x}}) ^\top {\bf A} ({\bf v}_1- \bar{{\bf x}})>0 \\ [6pt]
b=({\bf v}_2- \bar{{\bf x}}) ^\top {\bf A} ({\bf v}_2- \bar{{\bf x}})>0 \\ [6pt]
c=({\bf v}_1- \bar{{\bf x}}) ^\top {\bf A} ({\bf v}_2- \bar{{\bf x}}).
\end{array}
$$
By replacing (\ref{eq:optform}) in the objective function of (\ref{eq:optdir2}), we have that (\ref{eq:optdir2}) can be rewritten as
$$
\max_{\beta\in [0,1]} \min\left\{(-\beta^2+2 \beta) a-(1-\beta)^2 b+2(1-\beta)^2 c, -\beta^2 a +(1-\beta^2)b+ 2 \beta^2c       \right\}. 
$$
The optimal solution of this problem is
$$
\beta^\star=
\left\{
\begin{array}{ll}
0 & \mbox{if}\ b\leq c \\ [6pt]
1 &  \mbox{if}\ a\leq c \\ [6pt]
\frac{b-c}{a+b-2c} & \mbox{otherwise.}
\end{array}
\right.
$$
Then, the optimal value is:
$$
\left\{
\begin{array}{ll}
b & \mbox{if}\ b\leq c \\ [6pt]
a  &  \mbox{if}\ a\leq c \\ [6pt]
\frac{ab-c^2}{a+b-2c} & \mbox{otherwise.}
\end{array}
\right.
$$
We notice that $a,b>0$,
$$
a+b-2c=({\bf v}_1- {\bf v}_2) ^\top {\bf A} ({\bf v}_1- {\bf v}_2)>0,
$$
and, by Cauchy-Schwarz inequality:
$$
ab-2 c^2\geq 0,
$$
and equality holds if and only if $({\bf v}_1- \bar{{\bf x}})$ and $({\bf v}_2- \bar{{\bf x}})$ are linearly dependent. Thus, the optimal value of (\ref{eq:optdir2}) is always strictly positive unless 
the two vectors $({\bf v}_1- \bar{{\bf x}})$ and $({\bf v}_2- \bar{{\bf x}})$  lie along the same direction. 
More precisely, the optimal value is null only if the two vectors have the same direction but opposite sign. Indeed, let
$$
{\bf v}_1- \bar{{\bf x}}=\gamma({\bf v}_2- \bar{{\bf x}}).
$$
Then, we have $b=\gamma^2 a$ and $c=\gamma a$. If $\gamma$ is positive, then either $b\leq c$ (if $\gamma\leq 1$), or $a\leq c$ (if $\gamma\geq 1$) occurs, so that the optimal value is equal to $a$ or $b$ and is, thus, positive.
If $({\bf v}_1- \bar{{\bf x}})$ is not a negative multiple of $({\bf v}_2- \bar{{\bf x}})$, we are able to locally adjust $\bar{{\bf x}}$ along direction
$$
{\bf d}=\beta^\star ({\bf v}_1- \bar{{\bf x}}) +(1-\beta^\star)({\bf v}_2- \bar{{\bf x}}).
$$
\subsection{$P_1^\star(\lambda_{\Omega_{\bar{{\bf x}}}})$ is an infinite connected set}
\label{sec:x1=inf}
In this case we need to solve the following optimization problem
\begin{equation}
\label{eq:optdirinf}
\max_{{\bf d}\in \mathbb{R}^n} \min_{{\bf v}\in P_1^\star(\lambda_{\Omega_{\bar{{\bf x}}}})} -{\bf d} ^\top {\bf A} {\bf d}+ 2 {\bf d}^\top {\bf A} ({\bf v}- \bar{{\bf x}}).
\end{equation}
An improving direction exists if and only if the optimal value of this problem is strictly positive (note that the optimal value is always nonnegative since the inner minimization problem has optimal value 0 for ${\bf d}={\bf 0}$).
We first remark that the problem is convex. Indeed, for each fixed ${\bf v}$, we have a concave function with respect to ${\bf d}$, and the minimum of an infinite set of concave functions is itself a concave function (to be maximized, so that the problem is convex). 
The inner minimization problem can be solved in closed form. After removing the terms which do not depend on ${\bf v}$, the inner problem to be solved is
$$
 \min_{{\bf v}\in P_1^\star(\lambda_{\Omega_{\bar{{\bf x}}}})} 2 {\bf d}^\top {\bf A}  {\bf v}.
$$
According to Subcase 2.3.2 in Section \ref{sec:dual}, $P_1^\star(\lambda_{\Omega_{\bar{{\bf x}}}})$ can be written as in (\ref{eq:optimalset}) and the minimization problem can be reduced to the computation of the minimum of a linear function over the unit sphere:
$$
\min_{\boldsymbol{\xi}\in \mathbb{R}^q\ :\ \|\boldsymbol{\xi}\|^2=1} \bar{{\bf c}}({\bf d})^\top \boldsymbol{\xi},
$$
where $\bar{{\bf c}}({\bf d})$ is some linear function of ${\bf d}$ and $q\geq 2$ is the multiplicity of the minimum eigenvalue of the matrix ${\bf V}^\top ({\bf Q}+\lambda {\bf A}) {\bf V}$, corresponding to the Hessian of the objective function
of problem (\ref{eq:trustnull}).
The  optimal solution of this problem is 
$$
\boldsymbol{\xi}^\star=-\frac{\bar{{\bf c}}({\bf d})}{\|\bar{{\bf c}}({\bf d})\|},
$$
while the optimal value is $-\|\bar{{\bf c}}({\bf d})\|$.
\subsection{An algorithm for the refinement of the bound}
\begin{figure}[!tbh]
\centering
\includegraphics[width=.9\columnwidth]{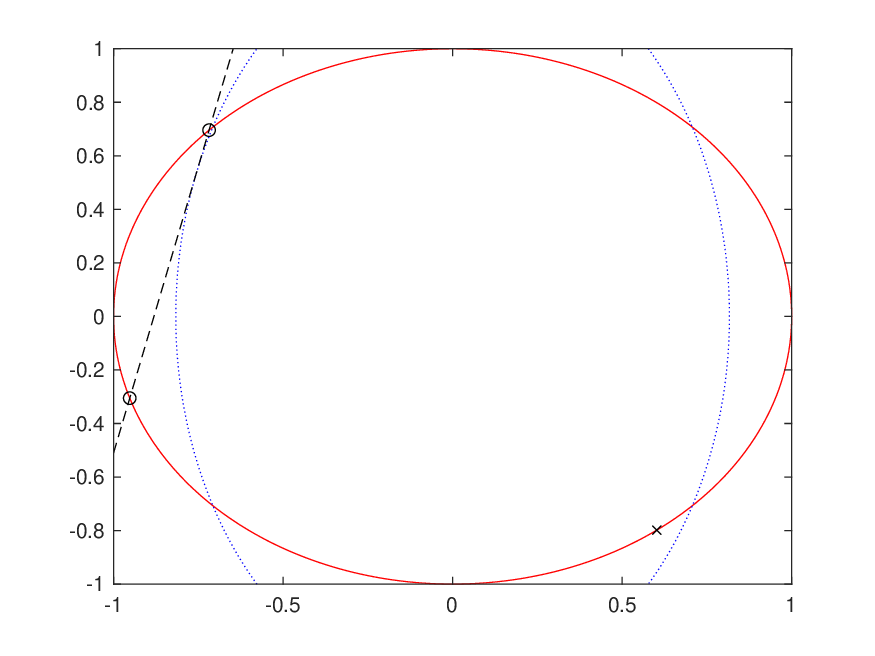}
\caption{Final linear cut after running  Algorithm \ref{algo:3}. Problem (\ref{eq:paramcdt1}) has got three optimal solutions, one in $int({H})$ and two outside ${H}$. The latter solutions are opposite to each other with respect to the final vector ${\bf z}$.}
\label{fig:twooptsolonecut}
\end{figure}

Let $\bar{\bf x}$ and $\lambda_{\mathbb{R}^n}$ be defined as in Section \ref{sec:boundimp}.
We propose Algorithm \ref{algo:3} for a bound based on successive local adjustments of the linear cut.
%Given the input point $\bar{\bf x}$ and the related $\bar{\lambda}$ value, CHIARIRE IN CHE SENSO RELATED
In line \ref{runalgo1first}, Algorithm \ref{algo:1} is run with input $X=\Omega_{\bar{\bf x}}$ and $\lambda_{\mathbb{R}^n}$. Note that with a slight abuse here we are assuming that the algorithm returns $\lambda_{\Omega_{\bar{\bf x}}}$ and the related points ${\bf z}_1$ and ${\bf z}_2$, while in practice close approximations of these quantities are returned, namely
$\lambda^{\max}$, ${\bf z}_1(\lambda^{\max})$ and ${\bf z}_2(\lambda^{\min})$.
In line \ref{initpoints}, ${\bf z}$ is initialized with the input point $\bar{\bf x}$ itself and the direction ${\bf d}^\star$, following the discussion in Section \ref{sec:x1=1},  is set equal to the difference between ${\bf z}_2(\lambda_{\Omega_{\bar{\bf x}}})$, the point outside $H$ returned by Algorithm \ref{algo:1}, and $\bar{\bf x}$. The outer while loop of the algorithm (lines \ref{outwhile1}-\ref{outwhile2}) is repeated until the bound is improved by at least a tolerance value $tol$. Inside this loop, in line \ref{instepsize}
the initial step size $\eta=1$ is set and a new incumbent ${\bf y}\in \partial H$ is computed. The inner while loop (lines \ref{innerwhile1}-\ref{innerwhile2}) computes the step size: until the optimal value of problem (\ref{eq:paramcdt1}) with $\bar{\bf x}={\bf y}$ and $\lambda=\lambda_{\Omega_{{{\bf z}}}}$, denoted by $opt$, is lower than the current lower bound $Lb$, 
%or there are points outside $H$ in the set of optimal solutions $P_{\Omega_{{\bf y}}}(\lambda_{\Omega_{{\bf z}}})$, 
we need to decrease the step size and recompute a new incumbent ${\bf y}$ (lines \ref{if1}-\ref{if2}). If the step size falls below a given tolerance value,
we exit the inner loop and also the outer one. 
Otherwise, we have identified a new valid incumbent 
and we set to 1 the exit flag $stop$ for the inner loop (line \ref{else}), so that, later on, a new linear inequality (\ref{eq:linearcut}) with $\bar{\bf x}={\bf y}$ will be computed. Then, at line \ref{algo2call} we run Algorithm \ref{algo:1} with input 
$X=\Omega_{{\bf y}}$ and $\lambda_{\Omega_{{\bf z}}}$. Finally, in line \ref{updpoint}, we update point ${\bf z}$ and the direction ${\bf d}^\star$.
We remark that at each iteration ${\bf z}_2(\lambda_{\Omega_{{\bf z}}})$ is {\em one} optimal solution of the current subproblem (\ref{eq:paramcdt1}) with $\lambda=\lambda_{\Omega_{{\bf z}}}$ lying outside $H$ and at which the linear cut of the subproblem is active, i.e.,  ${\bf z}_2(\lambda_{\Omega_{{\bf z}}})\in P_1^\star(\lambda_{\Omega_{{\bf z}}})$. As seen in Section \ref{sec:x1=1}, if $|P_1^\star(\lambda_{\Omega_{{\bf z}}})|=1$, i.e., ${\bf z}_2(\lambda_{\Omega_{{\bf z}}})$ is the unique optimal solution of the current subproblem (\ref{eq:paramcdt1}) with $\lambda=\lambda_{\Omega_{{\bf z}}}$  lying outside $H$, then, in view of Proposition \ref{prop:boundimprove}, 
the local adjustment employed in Algorithm \ref{algo:3} is guaranteed to improve the bound. However, as seen in Sections \ref{sec:x1=2} and \ref{sec:x1=inf}, if $P_1^\star(\lambda_{\Omega_{{\bf z}}})$ contains more than one point, than the proposed local adjustment
is not guaranteed to improve the bound. Sections  \ref{sec:x1=2} and \ref{sec:x1=inf} suggest how to define perturbing directions which still allow to improve the bound, in case they exist. However, as we will see through the computational experiments, Algorithm \ref{algo:3} turns  out to be time consuming and it is more convenient to improve the bound by adding a further linear cut, as we do in Section \ref{sec:two}, rather than further locally adjusting the current linear cut. In order to clarify this point,
we can make a comparison with Integer Linear Programming (ILP). In ILP problems, once a linear relaxation is solved, a valid cut removes {\em one} optimal solution of the relaxation. If the optimal solution is unique, then after the addition of the valid cut, the bound improves. But if the linear relaxation has got multiple solutions, then the valid cut is not guaranteed to remove all of them and, thus, the bound may not improve. It is possible to try to strengthen the valid cut in such a way that all optimal solutions of the linear relaxations are removed. But, more commonly, new linear cuts are added.
\begin{algorithm}
\caption{Bound improvement through a local adjustment of the linear cut.}
\label{algo:3}
\vspace{0.25cm}
%\hspace*{\algorithmicindent} 
\textbf{Input:} $\bar{{\bf x}}, \lambda_{\mathbb{R}^n}$
\vspace{0.25cm}
\begin{algorithmic}[1]
\STATE{Set $Lb_{old}=-\infty$}
%\STATE{Run Algorithm \ref{algo:1} with input $X=\Omega_{\bar{\bf x}}$ and $\bar{\lambda}$\label{runalgo1first}} 
\STATE{Let 
$[Lb,\lambda_{\Omega_{\bar{\bf x}}},{\bf z}_1(\lambda_{\Omega_{\bar{\bf x}}}), {\bf z}_2(\lambda_{\Omega_{\bar{\bf x}}})]= \textbf{DualLagrangian}(\Omega_{\bar{\bf x}},\lambda_{\mathbb{R}^n})$\label{runalgo1first}} 
%$[Lb, {\bf z}^\star, \lambda_{\Omega_{\bar{\bf x}}}]= \textbf{DualLagrangian}(\Omega_{\bar{\bf x}},\lambda_{\mathbb{R}^n})$\label{runalgo1first}} 
\STATE{Set ${\bf z}=\bar{\bf x}$ and ${\bf d}^\star={\bf z}_2(\lambda_{\Omega_{\bar{\bf x}}})-\bar{\bf x}$\label{initpoints}}
\WHILE{$Lb-Lb_{old}>tol$ \label{outwhile1}}% {\bf or} $\eta>tol_2$}
\STATE{Set $Lb_{old}=Lb$, $\eta=1$ and ${\bf y}= \Pi_{{\bf A}, {\bf a}}({\bf z}+{\bf d}^\star)  \in \partial H$ \label{instepsize}}
%}
\STATE{Set $stop=0$}
\WHILE{$stop=0$ and $\eta>\varepsilon$\label{innerwhile1}}
\STATE{Solve problem (\ref{eq:paramcdt1}) with $\bar{\bf x}={\bf y}$ and $\lambda=\lambda_{\Omega_{{\bf z}}}$ and let $opt$ be its optimal value}
%\IF{$opt < Lb$ {\tt or} $P_{\Omega_{{\bf y}}}(\lambda_{\Omega_{{\bf z}}})\setminus H\neq \emptyset$}
\IF{$opt < Lb$}
\STATE{Set $\eta=\eta/2$ \label{if1}} 
\STATE{Set ${\bf y}=\Pi_{{\bf A}, {\bf a}}({\bf z}+\eta {\bf d}^\star) \in \partial H$ \label{if2}}
%Set ${\bf x}_{\gamma}=\bar{{\bf x}}+\eta({\bf x}^\star-\bar{{\bf x}})-\gamma(2 {\bf A}\bar{{\bf x}}+{\bf a})$\; 
%Set $\bar{\gamma}=\min\{\gamma>0\ :\ {\bf x}_{\gamma}\in \partial {\cal E}\}$\; 
%Set ${\bf y}={\bf x}_{\bar{\gamma}}$\;
%\ELSE\STATE{Set $stop=1$\; Set ${\bf z}={\bf y}$ \label{else}}  
\ELSE\STATE{Set $stop=1$ \label{else}} 
\ENDIF  
\ENDWHILE\label{innerwhile2}
%Set $\lambda_{\min}=0$, $\lambda_{\max}=\bar{\lambda}$ \;
\IF{$stop=1$}
\STATE{Let $[Lb, \lambda_{\Omega_{{\bf y}}}, {\bf z}_1(\lambda_{\Omega_{{\bf y}}}), {\bf z}_2(\lambda_{\Omega_{{\bf y}}})]= \textbf{DualLagrangian}(\Omega_{{\bf y}},\lambda_{\Omega_{{\bf z}}})$\label{algo2call}}

\STATE{Set ${\bf z}={\bf y}$, ${\bf d}^\star={\bf z}_2(\lambda_{\Omega_{{\bf y}}})-{\bf z}$\;\label{updpoint}}
\ENDIF
\ENDWHILE \label{outwhile2}
\RETURN $Lb$
\end{algorithmic}
\end{algorithm}
$\ $\newline\newline\noindent
%In Figure \ref{fig:corrstep} we display the original point ${\bf v}$ for our example (denoted by $\circ$) (here $\eta=1$), the direction of the correction step (dotted line), the corrected point (denoted by $\times$), and, finally, the new linear cut (dashed line).
%\begin{figure}[!tbh]
%\centering
%\includegraphics[width=.9\columnwidth]{CorrectionStep.eps}
%\caption{Original point ${\bf v}$ ($\circ$), direction of the correction step (dotted line), corrected point ($\times$), and the new linear cut (dashed line).}
%\label{fig:corrstep}
%\end{figure}
Now we apply Algorithm \ref{algo:3} to our example.
\begin{example}
We have that ${\bf z}$ is initialized with $(-0.7901, 0.3565)$ and $Lb$ with $-4.0971$. 
During the execution of Algorithm \ref{algo:3}, ${\bf z}$ and $Lb$ are updated as indicated in Table \ref{tab:resiter}.
\begin{table}[ht]
    \centering
    \begin{tabular}{|c|c|c||}\hline
Iteration & ${\bf z}$ & $Lb$ \\\hline
1 &  $(-0.7204, 0.6658)$  &  $-4.0850$   \\ \hline 
2 &   $(-0.7742, 0.4493)$ &  $-4.0638$   \\ \hline 
3 &   $(-0.7481, 0.5665)$ &  $-4.0477$  \\ \hline 
4 &   $(-0.7607, 0.5136)$ &  $-4.0416$   \\ \hline 
5 &  $(-0.7556, 0.5361)$  &  $-4.0378$  \\ \hline 
6 &   $(-0.7571, 0.5296)$ & $-4.0364$    \\ \hline 
7 &  $(-0.7568, 0.5309)$  & $-4.0362$    \\ \hline 
\end{tabular}
    \caption{Iterations of Algorithm \ref{algo:3} over the example}
    \label{tab:resiter}
\end{table}
\end{example}
Interestingly, the best bound obtained in the example is exactly the one obtained for the same problem by the approach proposed in \cite{Burer13}, based on the addition of SOC-RLT constraints. 
Figure \ref{fig:twooptsolonecut} displays the situation at the last iteration of Algorithm \ref{algo:3}.
Problem (\ref{eq:paramcdt1}) has got three optimal solutions, one in $int(H)$ and two outside $H$. The two optimal solutions outside $H$ are opposite to each other with respect to the final vector ${\bf z}$, so that, as discussed in Section
\ref{sec:x1=2}, no further local adjustment is possible to improve the bound in this case.
%Thus, we conjecture that the bound in \cite{Burer13} is equal to the best possible bound which can be obtained by adding a single linear inequality, i.e., bound (\ref{eq:bestbound}).
\section{Bound improvement through the addition of a further linear cut}
\label{sec:two}
\begin{figure}[!tbh]
\centering
\includegraphics[width=.9\columnwidth]{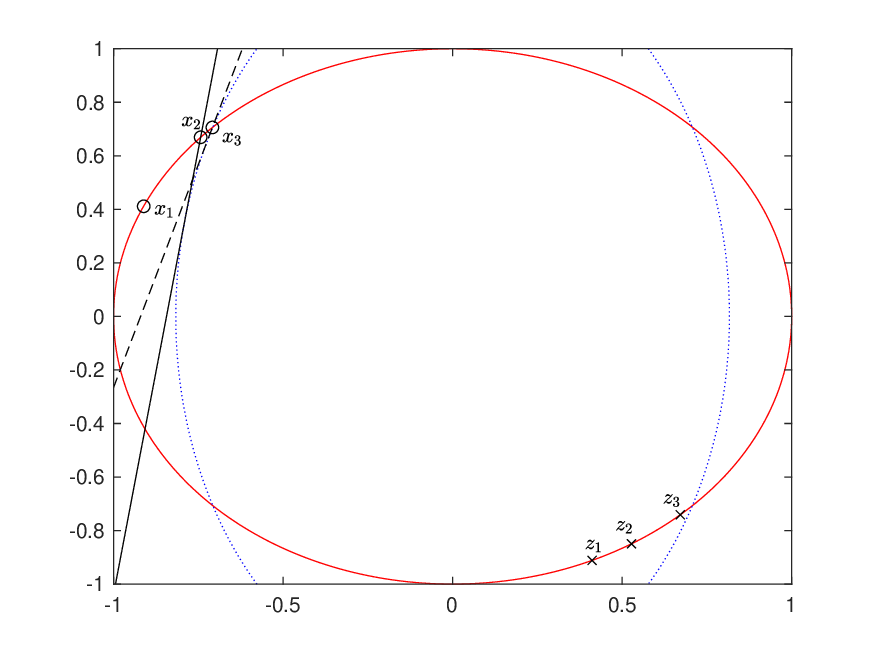}
\caption{Two linear cuts and the two optimal solutions outside ${H}$ ($x_3$) and in $int({H})$ ($z_3$), denoted by $\circ$ and $\times$, respectively.}
\label{fig:soltwolincut}
\end{figure}

\begin{figure}[!h]
\begin{subfigure}[b]{.6\textwidth}%{1.0\columnwidth}
\centering
\includegraphics[width=\columnwidth]{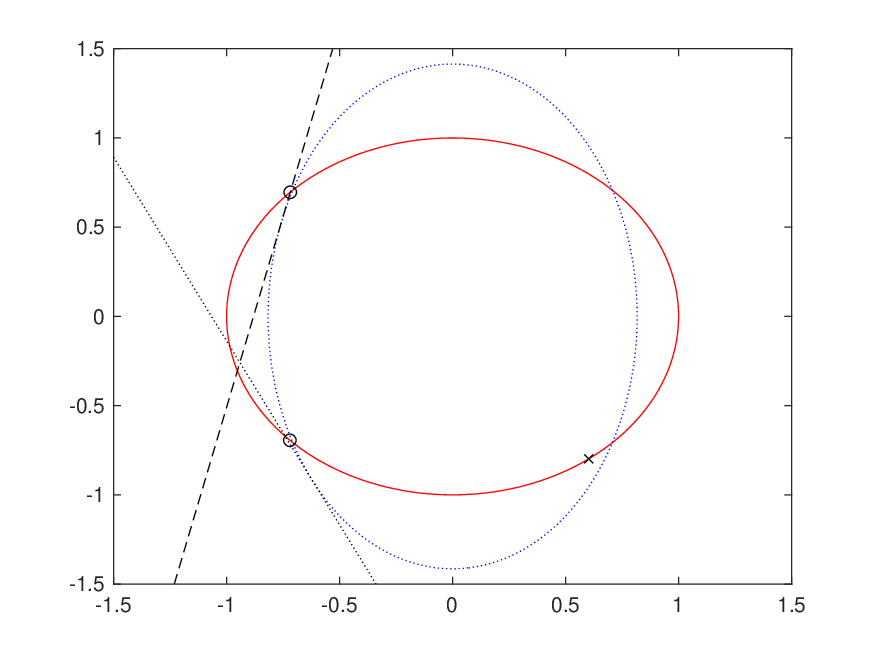}
\caption{Three optimal solutions, none with both linear cuts active.}
\label{fig:threeoptnoboth}
\end{subfigure}%
\newline
\begin{subfigure}[b]{.6\textwidth}%{0.5\columnwidth}
\centering
\includegraphics[width=\columnwidth]{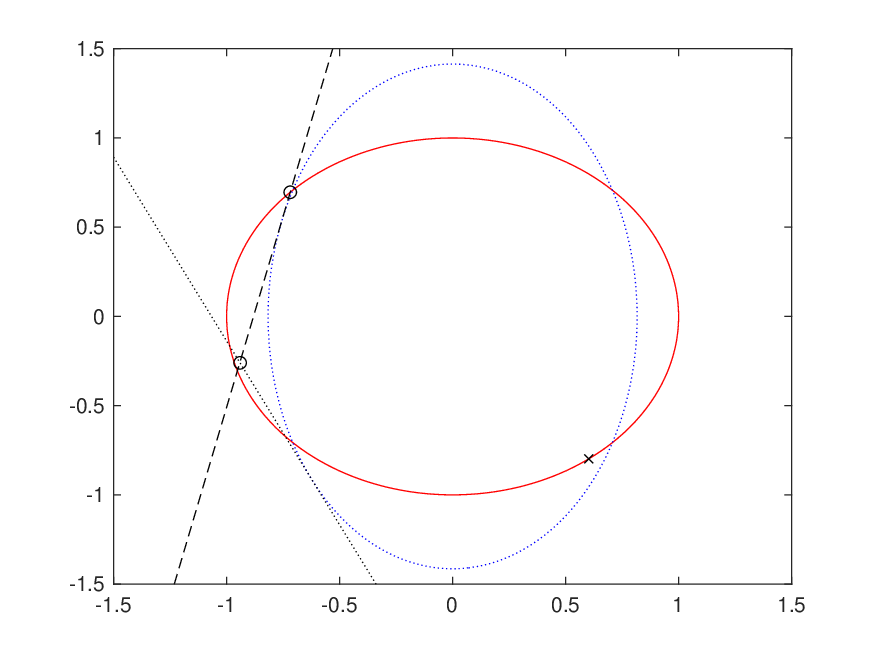}
\caption{Three optimal solutions, one with both linear cuts active.}
\label{fig:threeoptboth}
\end{subfigure}%
\newline
\begin{subfigure}[b]{.6\textwidth}%{0.5\columnwidth}
\centering
\includegraphics[width=\columnwidth]{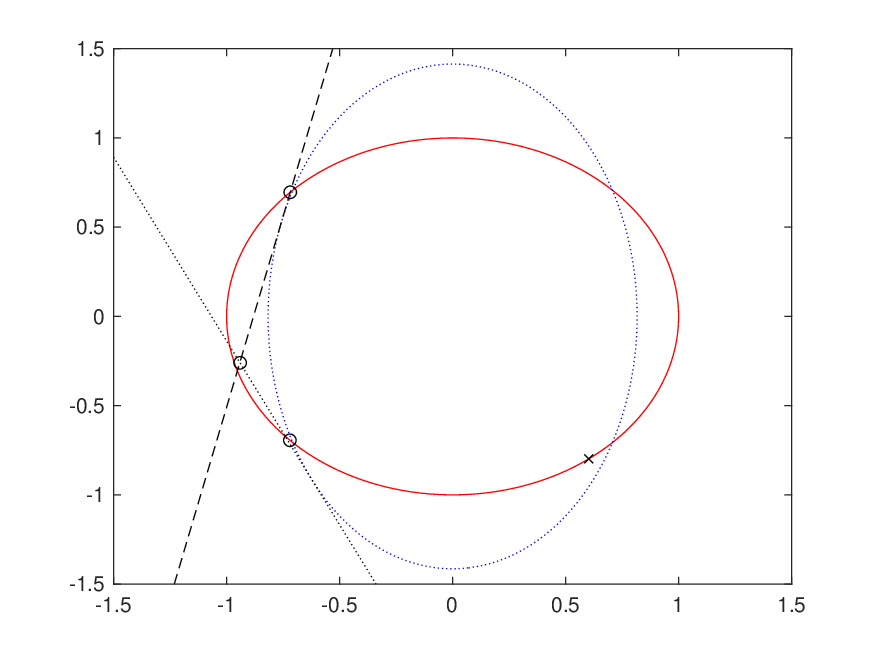}
\caption{Four optimal solutions.}
\label{fig:fouropt}
\end{subfigure}%
\end{figure}

\begin{figure}[!tbh]
\centering
\includegraphics[width=.9\columnwidth]{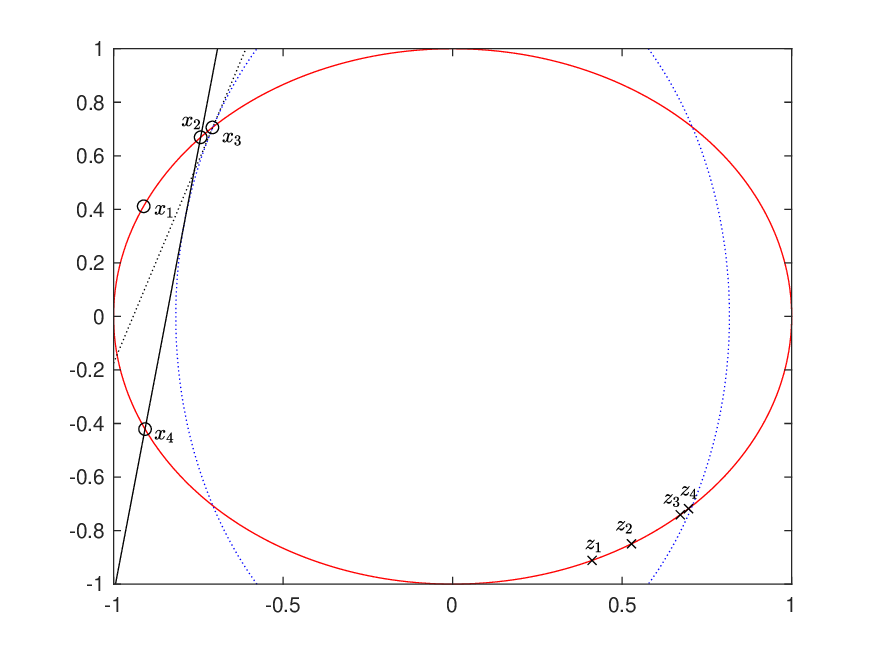}
\caption{Perturbation of the second linear cut and the two new optimal solutions outside ${H}$ ($x_4$) and in $int({H})$ ($z_4$), denoted by $\circ$ and $\times$, respectively.}
\label{fig:soltwofirstiter}
\end{figure}

\begin{figure}[!tbh]
\centering
\includegraphics[width=.9\columnwidth]{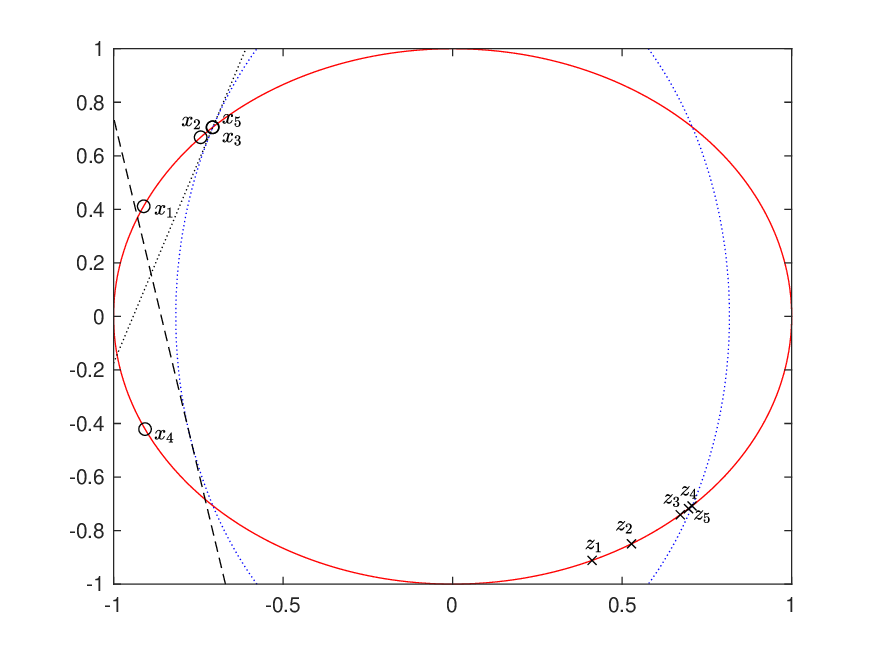}
\caption{Perturbation of the first linear cut and the two optimal solutions outside ${H}$ ($x_5$) and in $int({H})$ ($z_5$), denoted by $\circ$ and $\times$, respectively.}
\label{fig:soltwoseconditer}
\end{figure}

Another possible way to improve the bound is by adding a further linear cut to (\ref{eq:paramcdt1}). 
Let $\bar{\bf x}$ and $\lambda_{\mathbb{R}^n}$ be defined as in Section \ref{sec:boundimp}.
In line \ref{runalgo1first} of Algorithm \ref{algo:3}, we compute
$[Lb,\lambda_{\Omega_{\bar{\bf x}}},{\bf z}_1(\lambda_{\Omega_{\bar{\bf x}}}), {\bf z}_2(\lambda_{\Omega_{\bar{\bf x}}})]= \textbf{DualLagrangian}(\Omega_{\bar{\bf x}},\lambda_{\mathbb{R}^n})$,
%$[Lb,{\bf z}^\star,\lambda_{\Omega_{\bar{\bf x}}}]=\textbf{DualLagrangian}(\Omega_{\bar{\bf x}},\lambda_{\mathbb{R}^n})$ 
and, later on, we try to locally adjust $\bar{\bf x}$. Rather than doing that, we can add a further linear cut, cutting 
${\bf z}_2(\lambda_{\Omega_{\bar{\bf x}}})\not\in {H}$ away.
In particular, we add
the one obtained through the projection over $\partial {H}$ of ${\bf z}_2(\lambda_{\Omega_{\bar{\bf x}}})$. 
Let $\tilde{{\bf x}}= \Pi_{{\bf A}, {\bf a}}({\bf z}_2(\lambda_{\Omega_{\bar{\bf x}}}))\in \partial H$ be such projection. Then, we define the following problem
\begin{equation}
\label{eq:paramcdt2}
\begin{array}{ll}
\min_{{\bf x}} & {\bf x }^\top ({\bf Q}+\lambda {\bf A}) {\bf x}+({\bf q}+\lambda {\bf a})^\top {\bf x} -\lambda a_0 \\ [6pt]
 &  {\bf x}^\top {\bf x}\leq 1 \\ [6pt]
& (2{\bf A} \bar{{\bf x}}+{\bf a})^\top ({\bf x}-\bar{{\bf x}})\leq 0 \\ [6pt]
& (2{\bf A} \tilde{{\bf x}}+{\bf a})^\top ({\bf x}-\tilde{{\bf x}})\leq 0,
\end{array}
\end{equation}
which is equivalent to problem (\ref{eq:paramcdt}) where
$$
X=\Omega_{\bar{{\bf x}}}\cap \Omega_{\tilde{{\bf x}}}=\{{\bf x}\ :\ (2 {\bf A}\bar{{\bf x}}+{\bf a})^\top({\bf x}-\bar{{\bf x}})\leq 0,\ \ (2 {\bf A}\tilde{{\bf x}}+{\bf a})^\top({\bf x}-\tilde{{\bf x}})\leq 0\}\supset {H}.
$$ 
%Again Theorem \ref{lem:1} still holds true and $H_X^{\lambda},h_X^{\min}(\lambda),h_X^{\max}(\lambda)$ are always defined as
%in (\ref{eq:defdiffh}).
%However,
A convex reformulation as the one proposed in \cite{Burer13,Sturm03} for problem (\ref{eq:paramcdt1}) 
is not available in this case (unless the two linear inequalities do not intersect in the interior of the unit ball).
But in this case the alternative procedure discussed in Section \ref{sec:boundimp} turns out to be useful.
%Because of the addition of the new linear cut, we still have that for $\lambda$ values smaller but close to $\bar{\lambda}$ a unique optimal solution exists
%(QUI DISCUTERE IL CASO DI SOLUZIONI OTTIME MULTIPLE IN $F_1^\star(\lambda_{\min})$), belonging to the interior of ${H}$, so that also the linear constraints in (\ref{eq:paramcdt2}) are not active at it. Thus, such optimal solution must be a local and nonglobal minimizer of problem (\ref{eq:paramcdt}) (recall once again that the globally optimal solution of this problem always violates the second constraint in (\ref{eq:cdt}) for all $\lambda<\bar{\lambda}$).
As before, for each value $\lambda$ in the while loop of Algorithm \ref{algo:1} we can first check whether a local and nonglobal optimal solution of problem (\ref{eq:paramcdt}) with $X=\mathbb{R}^n$ exists, by exploiting the necessary and sufficient condition stated in \cite{WangXia20}.
%(recall once again that the globally optimal solution of this problem always violates the second constraint in (\ref{eq:cdt}) for all $\lambda\leq \bar{\lambda}^{\mathbb{R}^n}$).  QUESTO CHIARIRLO, FORSE NON E' CHIARO
%for each tested $\lambda$ value in the while loop we can first check whether a local and nonglobal optimal solution of problem (\ref{eq:paramcdt}) exists, by exploiting the necessary and sufficient condition stated in \cite{WangXia20}. 
If it exists, and belongs to $\Omega_{\bar{{\bf x}}}\cap \Omega_{\tilde{{\bf x}}}$, we denote it by ${\bf z}_1(\lambda)$. 
Next, we need to compute the optimal value
of (\ref{eq:paramcdt2}) when at least one of the two linear constraints is active, i.e., we need to solve the following problem
\begin{equation}
\label{eq:paramcdt3}
\begin{array}{ll}
\min_{{\bf x}} & {\bf x }^\top ({\bf Q}+\lambda {\bf A}) {\bf x}+({\bf q}+\lambda {\bf a})^\top {\bf x} -\lambda a_0 \\ [6pt]
 &  {\bf x}^\top {\bf x}\leq 1 \\ [6pt]
& (2{\bf A} \bar{{\bf x}}+{\bf a})^\top ({\bf x}-\bar{{\bf x}})\leq 0 \\ [6pt]
& (2{\bf A} \tilde{{\bf x}}+{\bf a})^\top ({\bf x}-\tilde{{\bf x}})\leq 0 \\ [6pt]
& \left[(2{\bf A} \bar{{\bf x}}+{\bf a})^\top ({\bf x}-\bar{{\bf x}})\right]\left[(2{\bf A} \tilde{{\bf x}}+{\bf a})^\top ({\bf x}-\tilde{{\bf x}})\right]=0.
\end{array}
\end{equation}
A convex reformulation of this problem has been proposed in \cite{Ye03}. Alternatively, one can solve two distinct problems, each imposing that one of the two linear inequalities is active.
Each of these problems can be converted into a trust region problem with an additional linear inequality, which can be solved in polynomial time through the already mentioned convex reformulation
proposed in \cite{Burer13,Sturm03}. 
Thus, we compute the set $P_1^\star(\lambda)\subseteq \partial \Omega_{\bar{\bf x}}\cap \Omega_{\tilde{\bf x}}$ of  optimal solutions of (\ref{eq:paramcdt2}) for which the first linear cut is active, and then 
the set $P_2^\star(\lambda)\subseteq \Omega_{\bar{\bf x}}\cap \partial \Omega_{\tilde{\bf x}}$ of  optimal solutions of (\ref{eq:paramcdt2})  for which the second linear cut is active.
Finally, the optimal values of these problems are compared with the value
of the local and nonglobal minimizer (if it exists) in order to identify the set $P_X(\lambda)$ of optimal solutions  of (\ref{eq:paramcdt2}). At this point we are able to compute  
$h_X^{\min}(\lambda),h_X^{\max}(\lambda)$ and update 
$\lambda^{\min}$ and $\lambda^{\max}$ accordingly.
%$\lambda^{\min}_{\Omega_{\tilde{\bf x}}\cap  \Omega_{\bar{\bf x}}}$ and $\lambda^{\max}_{\Omega_{\tilde{\bf x}}\cap  \Omega_{\bar{\bf x}}}$ accordingly. 
If for some $\lambda$ we have that ${\bf z}_1(\lambda)\in P_X(\lambda)$ and $P_X(\lambda)\cap [P_1^\star(\lambda)\cup P_2^\star(\lambda)]\neq \emptyset$, i.e., problem (\ref{eq:paramcdt2}) has an optimal solution in $int({H})$ and (at least) one optimal solution outside ${H}$, then $0\in [h^{\min}_X(\lambda),h^{\max}_X(\lambda)]$ and Algorithm~\ref{algo:1} stops.
%If the former is larger than the latter and the local and nonglobal minimizer belongs to ${H}$, then we set  $\lambda_{\max}=\lambda$, otherwise we set $\lambda_{\min}=\lambda$.
%Stated in another way, we still apply Algorithm \ref{algo:2}, with the only difference that now ${\bf z}_2(\lambda)$ is the solution of problem (\ref{eq:paramcdt3}), while ${\bf z}_1(\lambda)$ is still the local and nonglobal minimizer of the trust region problem (if it exists). 
%Algorithm \ref{algo:1}  will stop at some $\tilde{\lambda}<\bar{\lambda}$ at which again we have two optimal solutions, one of which violates the second constraint in (\ref{eq:cdt}).
We illustrate all this on Example \ref{ex:1}. 
\begin{example}
We add a second linear cut obtained through the projection over $\partial {H}$ of the optimal solution of problem (\ref{eq:paramcdt1}) with $\lambda_{\Omega_{\bar{\bf x}}}=0.726$ outside ${H}$. This leads to a further improvement  with 
 $\lambda_{\Omega_{\bar{\bf x}}\cap \Omega_{\tilde{\bf x}}}\approx 0.39$ and $p_{\Omega_{\bar{\bf x}}\cap \Omega_{\tilde{\bf x}}}(\lambda_{\Omega_{\bar{\bf x}}\cap \Omega_{\tilde{\bf x}}})\approx -4.005$, which almost closes the gap. 
In Figure \ref{fig:soltwolincut} we show the two linear cuts and the two new optimal solutions, one outside $H$ and one belonging to $int({H})$ ($x_3$ and $z_3$, respectively). Again, we also report the previous pairs of optimal solutions in order to show the progress.
\end{example}
Now, assume that the returned bound is not exact.
Also in this case $\bar{{\bf x}}$ and $\tilde{{\bf x}}$ can be locally adjusted. 
%For instance, in our example we could set $\tilde{{\bf x}}$ equal to the best local solution obtained by running two local searches for the original problem (\ref{eq:cdt}), one from
%the optimal solution of (\ref{eq:paramcdt1}) with the $\lambda$ value at the end of Algorithm \ref{algo:2} (namely, $\lambda=0.726$) belonging to the interior of ${H}$, and the other from the optimal solution of (\ref{eq:paramcdt1}) outside ${H}$.
%In this case, the bound is further reduced to
%$-4.0005$, very close to the true optimal value (-4), but there is still a positive (though small) gap. The positive gap of this small example suggests a further improvement. 
One can combine the techniques presented in Section \ref{sec:optimal} and in the current section, by using a technique similar to the one described in the former section to improve
the pair of points $\bar{{\bf x}}$ and $\tilde{{\bf x}}$. In particular, 
%let $\lambda_{\Omega_{\tilde{\bf x}}\cap  \Omega_{\bar{\bf x}}}$ be the limit value at which $\lambda^{\min}_{\Omega_{\tilde{\bf x}}\cap  \Omega_{\bar{\bf x}}}$ and $\lambda^{\max}_{\Omega_{\tilde{\bf x}}\cap  \Omega_{\bar{\bf x}}}$ converge when running Algorithm \ref{algo:1} with input $X=\Omega_{\bar{\bf x}}\cap \Omega_{\tilde{\bf x}}$. As previously remarked, 
at $\lambda_{\Omega_{\tilde{\bf x}}\cap  \Omega_{\bar{\bf x}}}$ we have one optimal solution of problem \ref{eq:paramcdt2} belonging to $int({H})$, namely the local and nonglobal optimal solution of problem (\ref{eq:paramcdt}) with $X=\mathbb{R}^n$,  and at least another one outside ${H}$. We denote the latter by ${\bf v}$ and we observe that at least one of the two linear cuts is active at this point, i.e., either ${\bf v}\in \partial \Omega_{\bar{\bf x}}$ or ${\bf v}\in \partial \Omega_{\tilde{\bf x}}$ (or both). Then, if only the first cut is active at ${\bf v}$, 
we update $\bar{\bf x}$ as follows 
$\bar{\bf x}'=\Pi_{{\bf A}, {\bf a}}(\bar{\bf x}+\eta ({\bf v}-\bar{\bf x}))$ for a sufficiently small $\eta$ value, while $\tilde{{\bf x}}'=\tilde{{\bf x}}$. 
If only the second cut is active, we update $\tilde{\bf x}$ as follows $\tilde{\bf x}'=\Pi_{{\bf A}, {\bf a}}(\tilde{\bf x}+ \eta({\bf v}-\tilde{\bf x}))$, for a sufficiently small $\eta$ value, while $\bar{\bf x}'=\bar{\bf x}$. 
Finally, if both are active we select one of the two cuts and perturb it. After the perturbation, we run again Algorithm \ref{algo:1} with input $X=\Omega_{\bar{{\bf x}}'}\cap \Omega_{\tilde{{\bf x}}'}$ and $\lambda_{\Omega_{\tilde{\bf x}}\cap  \Omega_{\bar{\bf x}}}$, and we repeat 
this procedure until there is a significant reduction of the bound. Note, however, that it might happen that no improvement is possible.  
%We illustrate different possible cases through Figures \ref{fig:threeoptnoboth}-\ref{fig:fouropt}.
%As usual, in these figures the point in the interior of ${H}$ is denoted by $\times$, while the others (outside ${H}$) are denoted by $\circ$.
In case $|P_1^\star(\lambda_{\Omega_{\tilde{\bf x}}\cap  \Omega_{\bar{\bf x}}})|=1$ and $P_2^\star(\lambda_{\Omega_{\tilde{\bf x}}\cap  \Omega_{\bar{\bf x}}})=\emptyset$ (similar for $|P_2^\star(\lambda_{\Omega_{\tilde{\bf x}}\cap  \Omega_{\bar{\bf x}}})|=1$ and $P_1^\star(\lambda_{\Omega_{\tilde{\bf x}}\cap  \Omega_{\bar{\bf x}}})=\emptyset$), then 
the proposed perturbation $\bar{\bf x}'=\Pi_{{\bf A}, {\bf a}}(\bar{\bf x}+\eta({\bf v}-\bar{\bf x}))$ for $\eta$ sufficiently small, allows to improve the bound. Indeed, in such cases the local adjustment is able to cut the unique solution outside ${H}$ away.
In order to  illustrate other different cases we employ Figures \ref{fig:threeoptnoboth}-\ref{fig:fouropt}.
As usual, in these figures the point in $int({H})$ is denoted by $\times$, while the others (outside ${H}$) are denoted by $\circ$.
If $|P_1^\star(\lambda_{\Omega_{\tilde{\bf x}}\cap  \Omega_{\bar{\bf x}}})|=|P_2^\star(\lambda_{\Omega_{\tilde{\bf x}}\cap  \Omega_{\bar{\bf x}}})|=1$ and $P_1^\star(\lambda_{\Omega_{\tilde{\bf x}}\cap  \Omega_{\bar{\bf x}}})\cap P_2^\star(\lambda_{\Omega_{\tilde{\bf x}}\cap  \Omega_{\bar{\bf x}}})=\emptyset$, 
(see Figure \ref{fig:threeoptnoboth}), or $|P_1^\star(\lambda_{\Omega_{\tilde{\bf x}}\cap  \Omega_{\bar{\bf x}}})|=2, |P_2^\star(\lambda_{\Omega_{\tilde{\bf x}}\cap  \Omega_{\bar{\bf x}}})|=1$ and $P_1^\star(\lambda_{\Omega_{\tilde{\bf x}}\cap  \Omega_{\bar{\bf x}}})\cap P_2^\star(\lambda_{\Omega_{\tilde{\bf x}}\cap  \Omega_{\bar{\bf x}}})\neq \emptyset$ (see Figure \ref{fig:threeoptboth}), then it is not possible to remove all the solutions outside ${H}$ by perturbing a single linear cut. Indeed, in both cases the perturbation of a single linear cut is able to remove just one of the two optimal solutions outside ${H}$. But it is possible to remove both by perturbing both linear cuts. Instead, Figure \ref{fig:fouropt} illustrates a case where
$|P_1^\star(\lambda_{\Omega_{\tilde{\bf x}}\cap  \Omega_{\bar{\bf x}}})|=|P_2^\star(\lambda_{\Omega_{\tilde{\bf x}}\cap  \Omega_{\bar{\bf x}}})|=2$ and $P_1^\star(\lambda_{\Omega_{\tilde{\bf x}}\cap  \Omega_{\bar{\bf x}}})\cap P_2^\star(\lambda_{\Omega_{\tilde{\bf x}}\cap  \Omega_{\bar{\bf x}}})\neq \emptyset$. In this case even the perturbation of both
linear cuts is unable to remove all three solutions outside ${H}$.
% This is indeed the situation occurring with the single instance for which the relative error is above $10^{-4}$. Note that such situation is the equivalent for the case of two linear cuts of the one discussed in Example \ref{ex:threeopt} for the case with a single linear cut.
The only way  to remove all three solutions outside ${H}$ is through the addition of a further linear cut, but, of course, this leads to a more complex problem with one trust region constraint and three linear inequalities.
%$\ $\newline\newline\noindent
\begin{example}
In our example, this refinement  is finally able to close the gap and return the exact optimal value $-4$.
In Figure \ref{fig:soltwofirstiter} we report the result of the first perturbation of the linear cuts. Since only the second linear cut is active at $x_3$, in this case the second linear cut is slightly perturbed and becomes equivalent to the tangent to ${H}$ at the optimal solution $(-\sqrt{2}/2,\sqrt{2}/2)$ of
the original problem (\ref{eq:cdt}). It is interesting to note that the new optimal solution outside ${H}$, indicated by $x_4$, lies in a different region with respect to the previous ones and is further from $\partial {H}$ with respect to
$x_2$ and $x_3$ (the reduction of $\lambda$ reduces the penalization of points outside ${H}$). Such solution is cut by the new linear inequality, obtained by a (not so small) perturbation of the first linear cut, displayed in Figure \ref{fig:soltwoseconditer}, together with the two new optimal solutions ($x_5$ and $z_5$), now corresponding to the two optimal solutions of problem (\ref{eq:cdt}).
\end{example}
%keeping $\tilde{{\bf x}}$ equal to a globally optimal solution and moving $\bar{{\bf x}}$ to point  $(-0.8005135, 0.278453)$ finally allows to close the gap.
\section{Computational Experiments}
\label{sec:comp}
In this section we report the computational results for the proposed bounds over the set of hard instances selected from the random ones generated in \cite{Burer13} and inspired by \cite{Martinez94}.
More precisely, in \cite{Burer13} 1000 random instances were generated for each size $n=5,10,20$. Some of these instances have been declared hard ones, namely those for which the bound obtained
by adding SOC-RLT constraints was not exact. In particular, these are 38 instances with $n=5$, 70 instances with $n=10$, and 104 instances with $n=20$. Such instances have been made available in {\tt GAMS}, {\tt AMPL}, and {\tt COCOUNT} formats
in \cite{Montahner18}. We tested our bounds on such instances.
All tests have been performed on an Intel Core i7 running at 1.8 GHz with 16GB of RAM. All bounds have been coded in {\tt MATLAB}.\newline\newline\noindent
We computed the following bounds:
\begin{itemize}
\item {\tt LbDual}, the dual Lagrangian bound computed through Algorithm \ref{algo:1} with input $X=\mathbb{R}^n$; 
\item {\tt LbOneCut}, the bound obtained by adding a single linear cut and computed through Algorithm \ref{algo:1} with input $X=\Omega_{\bar{\bf x}}$;
\item {\tt LbOneAdj}, the bound obtained by local adjustments of the added linear cut as indicated in Algorithm \ref{algo:3}; 
\item {\tt LbTwoCut}, the bound obtained by adding two linear cuts; 
%\item {\tt LbTwoCut-1}, the bound obtained by adding two linear cuts; 
%\item {\tt LbTwoCut-2}, the bound obtained by adding two linear cuts but the second one is obtained by taking the supporting hyperplane of the ellipsoid ${H}$ evaluated at the best local minimizer for (\ref{eq:cdt}) obtained by running two local searches, each one starting from one of the two optimal solutions obtained at the end of Algorithm \ref{algo:2}; 
\item {\tt LbTwoAdj} the bound obtained by adjusting the two linear cuts.
\end{itemize}
$\ $\newline\newline\noindent
According to what done in \cite{Anstreicher17,Burer13,Yang16}, an instance is considered to be solved when the
relative gap between the lower bound, say $LB$, and the upper bound, say $UB$, is not larger than $10^{-4}$, i.e.,
$$
\frac{UB-LB}{|UB|}\leq 10^{-4}.
$$
We set $UB$ equal to the lowest value obtained by running, after the addition of the first linear cut, two local searches for the original problem (\ref{eq:cdt}), one from
the optimal solution ${\bf z}_1(\lambda_{\Omega_{\bar{\bf x}}})\in int({H})$ of (\ref{eq:paramcdt1}) returned at the end of Algorithm \ref{algo:1}, and the other from an optimal solution of the same problem outside ${H}$.
In Tables \ref{tab:times5}-\ref{tab:times20} we report the average and maximum relative gaps for each bound, and the average and maximum computing times for $n=5,10,20$, respectively. 
Moreover, the average computing time for bound {\tt LbTwoAdj} is computed only over the instances (87 overall, as we will see) which are {\em not} solved by bound {\tt LbTwoCut}, while
for bounds {\tt LbTwoCut} and {\tt LbTwoAdj} the average gap is taken over the instances which were {\em not} solved by these bounds.
\newline\newline\noindent
We remark that the bound {\tt LbTwoCut} is computed by adding the first cut as in bound {\tt LbOneCut}, i.e., the supporting hyperplane at $\bar{\bf x}\in \partial {H}$, and then adding a further linear cut through the projection of an optimal solution outside ${H}$ obtained when computing 
bound {\tt LbOneCut}, i.e., point ${\bf z}_2(\lambda_{\min})$ returned by procedure ${\tt DualLagrangian}$ with input $X=\Omega_{\bar{\bf x}}$. We could as well choose the adjusted cut computed by bound {\tt LbOneAdj} as the first cut for bound {\tt LbTwoCut}, but we observed that with this choice no improvement over {\tt LbOneAdj} is obtained.
This is related to what already observed in Figure \ref{fig:twooptsolonecut}: bound {\tt LbOneAdj} cannot be improved any more when there are (at least) two optimal solutions outside ${H}$ (besides the one in $int({H}$)). Thus, the second cut
is able to remove one of such optimal solutions but not the other, so that the bound cannot be improved.
Similarly, for bound {\tt LbTwoAdj} the two initial cuts are the ones computed for bound {\tt LbTwoCut}.
\newline\newline\noindent
For what concerns the computing times, we observe that these are lower than those reported in \cite{Yang16} for the bound
obtained with the addition of SOC-RLT cuts (around 4s for an instance with $n=20$) and for the bound obtained by adding lifted-RLT cuts (around 92s for an instance with $n=20$).
They are also lower than those reported in \cite{Anstreicher17} for the bound obtained by adding KSOC cuts (up to 2s for $n=20$ instances). For the sake of correctness, we point out that the computing times reported in those papers have been obtained with different processors. However, such processors have comparable performance with respect to the one employed for the computational experiments in this paper.
In general, the proposed bounds are very cheap. Only for two instances with $n=20$, {\tt LbTwoAdj} required times above 1s (around 1.5s in both cases). Usually the computing times are (largely) below 1s. 
Both the dual Lagrangian bound and the bound obtained by a single linear cut are pretty cheap but with poorer performance in terms of relative gap. The bound obtained by Algorithm \ref{algo:3} 
with a local adjustment of the linear cut is better than the two previous ones in terms of gap but is also more expensive (although still cheap).
The bound {\tt LbTwoCut} offers a good combination between quality and cheap computing time. But a more careful choice of the two linear cuts, through a local adjustment, improves the quality without compromising the computing times.
This is confirmed by the results reported for ${\tt LbTwoAdj}$. Although this bound is more expensive than the others, the additional search for adjusted linear cuts further increases the quality of the bound.
In Table \ref{tab:solved} we report the number of solved instances for {\tt LbTwoCut}  and {\tt LbTwoAdj}. 
%Note that, according to what done in \cite{Anstreicher17,Yang16} an instance is considered to be solved when the
%relative difference between the lower bound, say $LB$, and the upper bound, say $UB$, is not larger than $10^{-4}$, i.e.,
%$$
%\frac{UB-LB}{|UB|}\leq 10^{-4}.
%$$
According to what reported in \cite{Anstreicher17}, the total number of unsolved instances out of the 212 hard instances is equal to:
133 for the bound proposed in \cite{Yang16}  (18 with $n=5$, 49 with $n=10$, and 66 with $n=20$);
85 for the bound proposed in \cite{Anstreicher17} (18 with $n=5$, 22 with $n=10$, and 45 with $n=20$);
56 by considering the best bound between the one in \cite{Yang16} and the one in \cite{Anstreicher17} (10 with $n=5$, 15 with $n=10$, and 31 with $n=20$).
For bound {\tt LbTwoCut} the total number of unsolved instances reduces to 87 (24, 29 and 34 for $n=5$, $n=10$, and $n=20$, respectively). 
Finally, for bound {\tt LbTwoAdj} we have the remarkable outcome that there is just one unsolved instance.
For the sake of correctness, we should warn that the value $UB$ in \cite{Anstreicher17,Yang16} is not computed by running two local searches as done in this paper.
It is instead computed from the final solution of the relaxed problem, so that it could be slightly worse and justify the larger number of unsolved instances.  
All the same, the quality of the proposed bounds appears to be quite good.

\begin{table}[ht]
    \centering
    \begin{tabular}{|l||c|c||c|c||}\hline
Bound & Average relative gap  (\%) &Max relative gap (\%) & Average time & Max time\\\hline
{\tt LbDual} & 0.90 \%   & 2.97 \% & 0.013  & 0.015 \\ \hline 
{\tt LbOneCut} & 0.31 \%     & 1.27 \%  & 0.035 & 0.040  \\ \hline 
{\tt LbOneAdj} & 0.13 \%      & 0.55 \%  &   0.266 & 0.388\\ \hline 
{\tt LbTwoCut} & 0.07 \%    & 0.21 \%   &  0.089 &  0.108 \\ \hline 
%{\tt LbTwoCut-1} & 0.07 \%    & 0.21 \%   &  0.089 &  0.108 \\ \hline 
%{\tt LbTwoCut-2} &  0.13 \%    & 0.26 \%   &  0.110 & 0.291 \\ \hline 
{\tt LbTwoAdj} & 0  \%    &  0 \%   & 0.146  & 0.281 \\ \hline 
\end{tabular}
    \caption{Average and maximum relative gaps and computing times (in seconds) for the instances with $n=5$}
    \label{tab:times5}
\end{table}

\begin{table}[ht]
    \centering
    \begin{tabular}{|l||c|c||c|c||}\hline
Bound & Average relative gap (\%)  &Max relative gap (\%) & Average time & Max time\\\hline
{\tt LbDual} &  0.41 \%   & 1.57 \%  & 0.014 & 0.022 \\ \hline 
{\tt LbOneCut} &  0.14 \%   & 0.81 \%  &  0.039& 0.057\\ \hline 
{\tt LbOneAdj} &   0.07 \%    & 0.48 \%  & 0.339  & 0.574 \\ \hline 
{\tt LbTwoCut} &  0.05 \%     & 0.24 \%  &   0.101 & 0.173\\ \hline 
%{\tt LbTwoCut-1} &  0.05 \%     & 0.24 \%  &   0.101 & 0.139\\ \hline 
%{\tt LbTwoCut-2} &   0.06 \%   & 0.18 \%  & 0.128 & 0.300 \\ \hline 
{\tt LbTwoAdj} & 0  \%    & 0 \%   & 0.197 & 0.670\\ \hline 
\end{tabular}
    \caption{Average and maximum relative gaps and computing times (in seconds) for the instances with $n=10$}
    \label{tab:times10}
\end{table}

\begin{table}[ht]
    \centering
    \begin{tabular}{|l||c|c||c|c||}\hline
Bound & Average relative gap (\%) &Max relative gap (\%) & Average time & Max time\\\hline
{\tt LbDual} & 0.20 \%   & 0.59 \%  &  0.019 & 0.027\\ \hline 
{\tt LbOneCut} &  0.08 \%    & 0.29 \%  &  0.057 & 0.079 \\ \hline 
{\tt LbOneAdj} &   0.05 \%    & 0.17 \%  &  0.539 & 0.926 \\ \hline 
{\tt LbTwoCut} &  0.03 \%     & 0.09 \%  &   0.148  & 0.199 \\ \hline 
%{\tt LbTwoCut-1} &  0.03 \%     & 0.09 \%  &   0.148  & 0.199 \\ \hline 
%{\tt LbTwoCut-2} & 0.08 \%    & 0.23 \%   &  0.196  & 0.481 \\ \hline 
{\tt LbTwoAdj} &  0.05 \%    &  0.05 \%   & 0.350 & 1.574\\ \hline 
\end{tabular}
    \caption{Average and maximum relative gaps and computing times (in seconds) for the instances with $n=20$}
    \label{tab:times20}
\end{table}

\begin{table}[ht]
    \centering
    \begin{tabular}{|l||c|c|c||}\hline
Bound &  $n=5$ (out of 38) &  $n=10$ (out of 70) &  $n=20$ (out of 104) \\\hline
%{\tt LbTwoCut-1} & 14  &41 & 70  \\ \hline 
{\tt LbTwoCut} & 14  &41 & 70  \\ \hline 
%{\tt LbTwoCut-2} & 35 & 63 & 89 \\ \hline 
{\tt LbTwoAdj} &  38     & 70  & 103   \\ \hline 
\end{tabular}
 %   \caption{Number of solved instances for the bounds {\tt LbTwoCut-1}, {\tt LbTwoCut-2} and {\tt LbTwoOpt}}.
\caption{Number of solved instances for the bounds {\tt LbTwoCut}  and {\tt LbTwoAdj}}.
    \label{tab:solved}
\end{table} 
\subsection{Investigating the hardest instance}
%VEDERE SE LASCIARE!!!!!!!!!!!!!!!!!
As a final experiment, we investigate the behaviour of bound {\tt LbTwoAdj} over the hardest instance with $n=20$, the one for which the relative error is above $10^{-4}$.
For this instance, at the last iteration we recorded the following objective function values, corresponding to values of local minimizers of problem (\ref{eq:paramcdt2}), which certainly include the global minimizer(s) of such problem:
\begin{itemize}
\item the value at the optimal solution of problem (\ref{eq:paramcdt2}) belonging to $int({H})$;
\item the value at a globally optimal solution of the trust region problem obtained by fixing in problem (\ref{eq:paramcdt2}) the first linear cut to an equality, in case such solution fulfills the second linear cut, or, alternatively, the value at the local and nonglobal solution of the same problem, in case such solution exists and fulfills the second linear cut (if the global minimizer does not fulfill the second linear cut and the local and nonglobal minimizer does not exist or does not fulfill the second linear cut, then the value is left undefined);
\item the same value as above but after fixing the second linear cut to an equality in problem (\ref{eq:paramcdt2});
\item the value at a globally optimal solution of the trust region problem obtained by fixing both cuts to equalities in problem (\ref{eq:paramcdt2}). 
\end{itemize}
Note that two of the four values must be equal. In particular, one of the two equal values is always the first one, attained 
in $int({H})$.
But for the hardest instance 
we observed that all four values are very close to each other and all of them are lower than the $UB$ value. Thus, it appears that for this instance a situation like the one
displayed in Figure \ref{fig:fouropt} occurs. In this case even the perturbation of both
linear cuts is unable to remove all the three solutions outside ${H}$. 

\section{Conclusions}
In this paper we discussed the CDT problem. 
First, we derived some theoretical results for a class of problems which includes the CDT problem as a special case. 
Then, from the theory developed for such class, we have re-derived a necessary and sufficient condition for the exactness of the Shor relaxation and of the equivalent dual Lagrangian bound for the CDT problem.
The condition is based on the existence of multiple solutions for a Lagrangian relaxation. Based on such condition, we proposed to strengthen the dual Lagrangian bound by adding one or two linear cuts.
These cuts are based on supporting hyperplanes of one of the two quadratic constraints and they are, thus, redundant for the original CDT problem (\ref{eq:cdt}). However, the cuts are not redundant for the Lagrangian relaxation and their
addition allows to improve the bound. We ran different computational experiments over the 212 hard test instances selected from the three thousand ones randomly generated in \cite{Burer13}, reporting gaps and computing times.
We have shown that the bounds are computationally cheap and are quite effective. In particular, one of them, based on the addition of two linear cuts, is able to solve all but one of the hard instances.
We have also investigated more in detail such hardest instance for which the bound is not exact (though quite close to the optimal value).
An interesting topic for future research could be that of establishing the relations between the bounds proposed in this work and those presented in the recent literature. Moreover, it would also be interesting to develop
procedures which are able to generate CDT instances for which the bound {\tt LbTwoAdj} is unable to return the optimal value.

%\section{Further developments}
%I would inspect in light of the current approach (addition of linear cuts to remove multiple solutions) the new bounds proposed in \cite{Yang16} and \cite{Yuan17} or, more generally, I would try to understand the relations between those results and the proposed approach.

\bibliographystyle{siamplain}

\begin{thebibliography}{9}
\bibitem{Adachi17}
S. Adachi, S. Iwata, Satoru, Y. Nakatsukasa, A.Takeda,
tSolving the trust-region subproblem by a generalized eigenvalue problem,
{\em SIAM Journal on Optimization}, 27(1), 269--291 (2017)
\bibitem{Ai09} W. Ai, S. Zhang, Strong duality for the CDT subproblem: A necessary and sufficient condition, 
{\em SIAM J. Optim.}, 19(4), 1735--1756 (2009)
\bibitem{Anstreicher17} K.M. Anstreicher, Kronecker product constraints with an application to the two-trust region subproblem, {\em SIAM J. Optim.}, 27(1), 368--378 (2017)
\bibitem{aubin2009set} J.P. Aubin, H. Frankowska, Set-Valued Analysis,
  {\em Modern Birkh{\"a}user Classics} (2009)
\bibitem{Barvinok93} A.I. Barvinok, Feasibility testing for systems of real quadratic equations, {\em Discrete Computational Geometry}, 10, 1--13 (1993)
%\bibitem{Bienstock13}   D.Bienstock, A. Michalka, Polynomial solvability of variants of the trust-region subproblem, {\em Proceedings of the 2014 Annual ACM-SIAM Symposium on Discrete Algorithms}
\bibitem{Beck06} A. Beck and Y. C. Eldar, Strong duality in nonconvex quadratic optimization with two
quadratic constraints, {\em SIAM J. Optim.}, 17, 844–860 (2006)
\bibitem{Berge1963}  C. Berge, Topological spaces : including a
  treatment of multi-valued functions, vector spaces and convexity, \emph{Dover Publications Inc.}. Mineola, New York (1963)
\bibitem{Bienstock16} D. Bienstock, A note on polynomial solvability of the CDT problem,
{\em SIAM Journal on Optimization}, 26, 488--498 (2016)
\bibitem{Bomze15} I.M. Bomze, M.L. Overton, Narrowing the difficulty gap
for the Celis-Dennis-Tapia problem, {\em Mathematical Programming}, 151, 459--476 (2015)  
\bibitem{Bomze18} I.M. Bomze, V. Jeyakumar, and G. Li, Extended trust-region problems with one or two balls: Exact copositive and Lagrangian relaxations, {\em Journal of Global Optimization}, 71, 551--569 (2018) 
\bibitem{Burer13} S.Burer, K.M. Anstreicher, Second-oder-cone constraints for extended trust-region subproblems, {\em SIAM J. Optim.}, 23(1), 432--451 (2013)
%\bibitem{BurerYang} S. Burer, B. Yang, The trust region subproblem with non-intersecting linear constraints, {\em Mathematical Programming}, 149, 253--264 (2015)
%\bibitem{Chen01} X.-D. Chen, Y.-X. Yuan, On maxima of dual function of the CDT subproblem, {\em Journal of Computational Mathematics}, 19(2) , 113--124 (2001)
\bibitem{Celis85} M.R. Celis, J.E. Dennis, R.A. Tapia, A trust region strategy for nonlinear equality constrained optimization.
In: Boggs, P.T., Byrd, R.H., Schnabel, R.B. (eds.) {\em Numerical Optimization 1984}, SIAM, Philadelphia
(1985)
\bibitem{Consolini17} L. Consolini and M.Locatelli, On the complexity of quadratic programming with two quadratic constraints,
{\em Mathematical Programming}, 164, 91--128 (2017)
\bibitem{Hil74} W. Hildenbrand, Core and Equilibria of a Large
    Economy. (PSME-5),\emph{Princeton University Press} (1974).
\bibitem{ConAl} J.B. Hiriart-Urruty, C. Lemarechal, Convex Analysis
  and Minimization Algorithms I - Fundamentals \emph{Springer} (1993)
\bibitem{Martinez94} J.M. Martinez, Local minimizers of quadratic functions on Euclidean balls and spheres,{\em SIAM J. Optim.}, 4(1), 159--176 (1994)
\bibitem{Montahner18} T. Montahner, A. Neumaier, F. Domes, A computational study of global optimization solvers on two trust region subproblems, {\em Journal of Global Optimization}, 71, 915--934 (2018)
\bibitem{More83} J.J. Mor\'e and D. C. Sorensen, Computing a trust region step,
{\em SIAM Journal on Scientific and Statistical Computing}, 4, 553--572 (1983) 
\bibitem{Sorensen82} D. C. Sorensen, Newton's  method  with  a  model  trust  region  modification, {\em SIAM J. Numer.Anal.}, 19, 409--426 (1982)
\bibitem{Sakaue16} S. Sakaue,  Y. Nakatsukasa, A.Takeda, A. and S. Iwata, Solving generalized CDT problems via two-parameter eigenvalues,
{\em SIAM Journal on Optimization}, 26, 1669--1694 (2016)
\bibitem{Sturm03} J.F. Sturm, S. Zhang, On cones of nonnegative quadratic functions, {\em Mathematics of Operations Research}, 28(2), 246--267 (2003)
\bibitem{WangXia20} J. Wang, Y. Xia, Closing the gap between necessary and sufficient conditions for local nonglobal minimizer of trust region subproblems, {\em SIAM J. Optim.}, 30(3), 1980--1995 (2020)
\bibitem{Yang16} B. Yang, S. Burer, A two-varable approach to the two-trust region subproblem, {\em SIAM J. Optim.}, 26(1), 661--680 (2016)
\bibitem{Ye03} Y. Ye, S. Zhang, New results on quadratic minimization,  {\em SIAM J. Optim.}, 14(1), 245--267 (2003)
\bibitem{Yuan17}  J. Yuan, M. Wang, W. Ai, T. Shuai, New results on narrowing the duality gap on the extended Celis-Dennis-Tapia problem, {\em SIAM J. Optim.}, 27(2), 890--909 (2017)
\end{thebibliography}

\end{document}